\documentclass[12pt,fullpage]{article}

\usepackage{amsfonts,graphicx,amsmath,xr,epsfig,comment}
\input epsf

\def\Z{\mathbb{Z}}
\def\R{\mathbb{R}}
\def\N{\mathbb{N}}
\def\epsilon{\varepsilon}
\def\trait (#1) (#2) (#3){\vrule width #1pt height #2pt depth #3pt}
\def\fin{\hfill\trait (0.1) (5) (0) \trait (5) (0.1) (0) \kern-5pt \trait (5) (5) (-4.9) \trait (0.1) (5) (0)}

\newcommand{\SE}{\setcounter{equation}{0} \section}
\newcommand{\be}{\begin{equation}}
\newcommand{\ee}{\end{equation}}
\newcommand{\baa}{\begin{array}}
\newcommand{\eaa}{\end{array}}
\newcommand{\ba}{\begin{eqnarray}}
\newcommand{\ea}{\end{eqnarray}}

\newtheorem{theo}{\bf Theorem}[section]
\newtheorem{lem}[theo]{\bf Lemma}

\newtheorem{defi}[theo]{\bf Definition}
\newtheorem{rem}[theo]{\bf Remark}

\begin{document}
\date{}
\title{\bf{Spreading speeds for some reaction-diffusion equations with general initial conditions}}
\author{Fran\c cois Hamel$^{\hbox{\small{ a, b}}}$\thanks{The first author is indebted to the Alexander von~Humboldt Foundation for its support. The two authors are also supported by the ANR project PREFERED.}  $\ \!$ and Yannick Sire$^{\hbox{\small{ a, c}}}$\\
\\
\footnotesize{$^{\hbox{a }}$ Aix-Marseille Universit\'e, LATP, Facult\'e des Sciences et Techniques}\\
\footnotesize{Avenue Escadrille Normandie-Niemen, F-13397 Marseille Cedex 20, France}\\
\footnotesize{$^{\hbox{b }}$ Helmholtz Zentrum M\"unchen, Institut f\"ur Biomathematik und Biometrie}\\
\footnotesize{Ingolst\"adter Landstrasse 1, D-85764 Neuherberg, Germany}\\
\footnotesize{$^{\hbox{c }}$ Laboratoire Poncelet,UMI $2615$, $119002$, Bolshoy Vlasyevskiy Pereulok 11, Moscow, Russia}}
\maketitle

\begin{abstract} 
This paper is devoted to the study of some qualitative and quantitative aspects of nonlinear propagation phenomena in diffusive media. More precisely, we consider the case a reaction-diffusion equation in a periodic medium  with ignition-type nonlinearity, the heterogeneity being on the nonlinearity, the operator and the domain. Contrary to previous works, we study the asymptotic spreading properties of the solutions of the Cauchy problem with general initial conditions which satisfy very mild assumptions at infinity. We introduce several concepts generalizing the notion of spreading speed and we give a complete characterization of it when the initial condition is asymptotically oscillatory at infinity. Furthermore we construct, even in the homogeneous one-dimensional case, a class of initial conditions for which highly nontrivial dynamics can be exhibited.  
\end{abstract}

%%%%%%%%%%%%%%%%%%%%%%%%%%%%%%%%%%%%%%%%
%%%%%%%%%%%%%%%%%%%%%%%%%%%%%%%%%%%%%%%%

\SE{Introduction}\label{intro}

We consider reaction-diffusion-advection equations of the type
\begin{equation}\label{per}
\left \{
\begin{array}[c]{rcll}
u_t-\nabla\cdot(A(z)\nabla u)+q(z)\cdot\nabla u & = & f(z,u), & \mbox{$z\in\overline{\Omega}$},\vspace{5pt}\\
\nu A\nabla u & = & 0, & \mbox{$z\in\partial\Omega$}
\end{array}\right.
\end{equation}
in an unbounded domain (connected and open) $\Omega\subset\R^N$ which is of class $C^{2,\alpha}$ for some $\alpha>0$. We denote by $\nu$ the outward unit normal on $\partial\Omega$. For any two vectors $\xi=(\xi_1,\ldots,\xi_N)$ and $\xi'=(\xi'_1,\ldots,\xi'_N)$ in $\R^N$ and any $N\times N$ matrix $B=(B_{ij})_{1\le i,j\leq N}$ with real entries, we write
$$\xi B\xi'=\sum_{1\le i,j\le N}\xi_iB_{ij}\xi'_j.$$

The coefficients of \eqref{per} are not assumed to be homogeneous in general, as well as the underlying domain $\Omega$. Instead, we just assume that there exist two real numbers $L>0$ and $R\ge 0$ such that
\be\label{omega}\left\{
\begin{array}[c]{c}
\forall\ z=(x,y)\in\Omega,\quad |y|\le R,\vspace{5pt}\\
\forall\ k\in L\mathbb{Z}\times\{0\}^{N-1},\quad\Omega=\Omega+k,\end{array}\right.
\ee
where 
$$x=x_1,\quad y=(x_{2},\cdots,x_N),\quad z=(x,y)$$
and $|\!\cdot\!|$ denotes the euclidean norm. The domain $\Omega$ is then an infinite cylinder which is unbounded in the direction $x$, its boundary $\partial\Omega$ may be straight or undulating, and $\Omega$ may also contain periodic perforations. Let $C$ be the periodicity cell defined by
$$C=\{z=(x,y)\in\Omega,\ x\in(0,L)\}.$$
Throughout the paper, we assume that the matrix field $z\mapsto A(z)=(A_{ij}(z))_{1\le i,j\le N}$ is of class $C^{1,\alpha}(\overline{\Omega})$, symmetric and uniformly elliptic, that the vector field $z\mapsto q(z)=(q_i(z))_{1\le i\le N}$ is of class $C^{0,\alpha}(\overline{\Omega})$, divergence-free ($\nabla \cdot q=0$ in $\overline{\Omega}$) and tangent to the boundary of $\Omega$ ($q\cdot\nu=0$ on $\partial\Omega$), and that the nonlinearity $f:(z,u)\ (\in\overline{\Omega}\times\R)\ \mapsto f(z,u)$ is continuous, of class $C^{0,\alpha}$ with respect to $z$ locally uniformly in $u\in\R$ and we assume that the restriction of $f$ to $\overline{\Omega}\times[0,1]$ is of class $C^1$ with respect to $u$. All functions $A_{ij}$, $q_i$ and $f(\cdot,u)$ (for all $u\in\R$) are assumed to be periodic, in the sense that they satisfy
$$w(x+k,y)=w(x,y)\ \hbox{ for all }(x,y)\in\Omega\hbox{ and }k\in L\mathbb{Z}.$$
The vector field $q$ is normalized in such a way that
$$\int_Cq(z)dz=0.$$
The nonlinearity $f$ is also assumed to be of combustion type: there exists $\theta \in (0,1)$ such that for every $z \in\overline{\Omega}$, 
\begin{equation}\label{comb}\left\{\begin{array}{l}
f(z,\cdot) \equiv 0 \,\,\,\mbox{on}\,\, (-\infty, \theta] \cup
[1,+\infty),\vspace{5pt}\\
f(z,\cdot)>0\,\,\,\,\mbox{on}\,\,(\theta,1),\,\,\,\,\displaystyle{\frac{\partial f}{\partial u}}(z,1^-)=-\displaystyle{\mathop{\lim}_{s\to 0^+}}\displaystyle{\frac{f(z,1-s)}{s}}<0.\end{array}\right.  
\end{equation}

Under the previous structural assumptions on the domain and the nonlinearity, we study the Cauchy problem 
\begin{equation}\label{problem}
\left \{
\begin{array}[c]{rcll}
u_t-\nabla\cdot(A(z)\nabla u)+q(z)\cdot\nabla u & = & f(z,u), & t>0,\ \mbox{$z\in\overline{\Omega}$},\vspace{5pt}\\
\nu A\nabla u & = & 0, & t>0,\ \mbox{$z\in\partial\Omega$},\vspace{5pt}\\
u(0,z) & = & u_0(z), & \mbox{$z\in \Omega$},
\end{array} \right . 
\end{equation}
where the initial value $u_0: \Omega \rightarrow [0,1]$ is uniformly continuous and satisfies the following mild conditions at infinity 
$$\displaystyle{\mathop{\lim}_{A \rightarrow +\infty}}\Big(\displaystyle{\mathop{\sup}_{z=(x,y) \in \Omega, \,x \geq A}}u_0(z)\Big)<\theta\ \hbox{ and }\ \displaystyle{\mathop{\lim}_{A \rightarrow -\infty}}\Big(\displaystyle{\mathop{\inf}_{z=(x,y) \in \Omega, \,x \leq A}}\ u_0(z)\Big)>\theta.$$
For sake of conciseness, we will denote the previous limits as
\begin{equation}\label{lim12}
\displaystyle{\mathop{\limsup}_{ x \rightarrow +\infty}}\ u_0(z) < \theta\ \hbox{ and }\ \displaystyle{\mathop{\liminf}_{ x \rightarrow -\infty}}\ u_0(z) > \theta.
\end{equation} 
The assumption of uniform continuity for $u_0$ is just made to ensure the solvability of the Cauchy problem. Notice also that, since $u_0$ satisfies $0\le u_0\le 1$ in $\Omega$ and is not identically equal to $0$ or $1$ because of (\ref{lim12}), the solution $u$ of (\ref{problem}) satisfies
\be\label{u01}
0<u(t,z)<1\ \hbox{ for all }t>0\hbox{ and }z\in\overline{\Omega}
\ee
from the strong parabolic maximum principle and Hopf lemma.

The main assumption (\ref{lim12}) means, roughly speaking, that the initial condition $u_0$ is front-like in the direction $x$, uniformly with respect to the orthogonal variables $y$. But it is important to notice that we do not assume that $u_0$ converges to some constants as $x\to\pm\infty$. The goal of this paper is to study propagation phenomena for the solutions $u$ of \eqref{problem} when the initial conditions $u_0$ just satisfy \eqref{lim12}. We shall see that these very weak assumptions at initial time give rise to a large variety of asymptotic spreading properties and possibly complex large-time behaviour. To this end, we first define the following two quantities, which shall stand for minimal and maximal asymptotic spreading speeds:

\begin{defi}\label{def1}
Let $u_0$ be as before. We define the {\rm lower spreading speed} $c_*(u_0)$ associated to \eqref{problem} as 
$$c_*(u_0)=\sup \mathcal E_*(u_0)$$
where 
$$\mathcal E_*(u_0)= \left \{ c\in \R \ |\ \forall\, c'<c,\  
\lim_{t \rightarrow +\infty}\Big(\inf_{ z \in \overline{\Omega},\ \!x \leq c't } u(t,z)\Big)=1 \right \}. $$
We also define the {\rm upper spreading speed} $c^*(u_0)$ associated to \eqref{problem} as 
$$c^*(u_0)=\inf \mathcal E^*(u_0)$$
where 
$$\mathcal E^*(u_0)= \left \{ c\in \R \ |\ \forall\, c'>c,\
\limsup_{t \rightarrow +\infty}\Big(\sup_{ z \in \overline{\Omega},\ \!x \geq c't } u(t,z)\Big)< 1 \right \}.$$  
\end{defi}

Qualitatively, the previous definitions show that an observer who moves at speed $c$ in direction $x$ will see for large times the steady state $1$ if $c <c _*(u_0)$ and will be away from~$1$ if $c > c^*(u_0)$. It follows in particular from Definition~\ref{def1} and (\ref{u01}) that, for all $A\in\R$,
$$\lim_{t \rightarrow +\infty}\Big(\sup_{(x+s,y) \in \overline{\Omega},\ \!x\le A,\ \!s\leq ct} |u(t,x+s ,y)-1|\Big)=0\,\,\,\,\,\mbox{if }c <c_*(u_0)$$
and
$$\limsup_{t \rightarrow +\infty}\Big(\sup_{ (x+s ,y) \in \overline{\Omega},\ \!x\ge A,\ \!s\geq ct} u(t,x+s ,y )\Big)<1\,\,\,\,\mbox{if }c >c^*(u_0).$$
Notice that by definition, there always holds 
$$c_*(u_0) \leq c^*(u_0).$$
Furthermore, if $c_*(u_0) \in \R $, resp. $c^*(u_0) \in \R$ --we shall see in Theorem~\ref{mainTh} that this is automatically true due to (\ref{lim12})-- then $c_*(u_0)=\max \mathcal E_*(u_0)$, resp. $c^*(u_0)=\min \mathcal E^*(u_0)$. However, this does not mean in general that 
$$\lim_{t \rightarrow +\infty}\Big(\inf_{ z \in \overline{\Omega},\ \!x \leq c_*(u_0)t } u(t,z)\Big)=1$$
or 
$$\limsup_{t \rightarrow +\infty}\Big(\sup_{ z \in \overline{\Omega},\ \!x \geq c^*(u_0)t } u(t,z)\Big)<1.$$

This paper is devoted to some characterizations of the lower and upper spreading speeds~$c_*(u_0)$ and $c^*(u_0)$ given in Definition~\ref{def1}, when $u_0$ satisfies the above conditions~\eqref{lim12}. We will derive some estimates for these spreading speeds and provide an example for which $c^*(u_0) \neq c_*(u_0)$, even in the homogeneous case.

One of the key points to understand propagation phenomena for the Cauchy pro\-blem~\eqref{problem} is based on the existence of a family of pulsating travelling fronts for \eqref{per}. In particular, we shall relate the spreading speeds $c_*(u_0)$ and $c^*(u_0)$ to various speeds of fronts connecting two stationary states of the equation. Given any real number $\gamma\in(-\infty,\theta)$, a pulsating front connecting $\gamma$ to $1$ and travelling to the right with effective speed $c\neq 0$ is a special time-global solution $u:\R\times\overline{\Omega}\to(\gamma,1)$ of~\eqref{per} satisfying the periodicity condition
\be\label{ptf}
\forall \, k \in L \mathbb Z ,\ \forall\, (t,x,y) \in \R \times \overline{\Omega},\quad u\Big(t-\frac{k}{c},x,y\Big)=u(t,x+k,y)
\ee
and the additional asymptotic conditions 
\be\label{limits}
\lim_{x  \rightarrow +\infty} u(t,z)=\gamma\ \hbox{ and }\ \lim_{x  \rightarrow -\infty} u(t,z)=1.
\ee
The previous limits (\ref{limits}) are taken locally in time and uniformly in $y$. Another way to describe a pulsating front is to use a hull function $\varphi:\R\times\overline{\Omega}\mapsto(\gamma,1)$ and write $u$ as
$$u(t,z)=\varphi(x-ct,z)$$
where the function $z (\in\overline{\Omega}) \mapsto\varphi(s,z)$ is periodic in $\overline{\Omega}$ for each $s\in\R$, and
$$\varphi(+\infty,\cdot)=\gamma,\ \ \varphi(-\infty,\cdot)=1\ \hbox{ uniformly in }\overline{\Omega}.$$
The existence and properties of pulsating travelling fronts have been obtained in \cite{x2,x3} for the case of the whole space $\R^N$ and in \cite{bh,bhbook} in the general periodic framework and with general combustion-type nonlinearities, covering the situation of the present paper. We sum up the result in the following theorem 

\begin{theo}\label{thFront} {\rm{\cite{bh,bhbook}}}
Let the nonlinearity $f$ be of the combustion type $(\ref{comb})$. For any $\gamma\in(-\infty,\theta)$, there exists a unique  speed $c=c_\gamma$, which is positive, such that problem~\eqref{per} has a pulsating travelling front solution $u_{\gamma}$ satisfying $(\ref{ptf})$ and $(\ref{limits})$. Furthermore, the function $u_{\gamma}$ is unique up to shifts in time and the map $\gamma \mapsto c_\gamma$ is continuous and increasing. 
\end{theo}

Under assumption (\ref{comb}), a great attention has been to be devoted to the properties of fronts in the homogeneous one-dimensional version of (\ref{per}), and then in straight infinite cylinders (see e.g. \cite{bn}). Of particular interest are the stability of these fronts and the convergence to the fronts of the solutions of Cauchy problems of the type \eqref{problem} when the initial condition $u_0$ is in some sense close to a given front and has the same (constant) limit as it when $x\to+\infty$ \cite{blr,k1,r1,r2,r3}. Initial conditions with compact support have also been considered. Under some conditions, that is if they are above and away from $\theta$ on a sufficiently large set, then they develop into a pair of diverging fronts \cite{k2,r3,z1}. However, in the general periodic setting, the question of the global stability of the travelling fronts still remains open, even for initial conditions having the same constant limit $\gamma$ as a given front $u_{\gamma}$ when $x\to+\infty$. As a matter of fact, the present paper is at least twofold: firstly, we show the convergence in speed for a more general class of asymptotically periodic (when $x\to+\infty$) initial conditions, and secondly we prove that such convergence does not hold in general, even in the homogeneous case, when the initial conditions just satisfy (\ref{lim12}).

Let us mention here that other types of nonlinearities have also been considered in the literature. For instance, some existence and stability results of fronts with bistable reaction terms are known, but they are mainly concerned with homogeneous or close-to-homogeneous media, or with media which are invariant in the direction of propagation \cite{fm,lx,r3,x1}. One of the most famous results in this spirit is the following one: in the homogeneous setting with bistable reaction-terms $f:[0,1]\to\R$ satisfying
$$f(0)=f(\theta)=f(1)=0,\ f<0\hbox{ on }(0,\theta),\ f>0\hbox{ on }(\theta,1),\ f'(0)<0,\ f'(1)<0$$
for some $\theta\in(0,1)$, front-like initial conditions satisfying (\ref{lim12}) are known to converge to the unique front connecting the two stable zeroes $0$ and $1$ of $f$ \cite{fm}. This is due to the strong attractivity of these two stable states. As will be seen, in the combustion case (\ref{comb}) considered in the present paper, new interesting and more complex phenomena shall occur, due to the existence of a continuum of stationary states (below $\theta$). Lastly, for monostable or particular Kolmogorov-Petrovski-Piskunov \cite{kpp} type nonlinearities, existence and qualitative properties of pulsating travelling fronts in periodic media have been established in \cite{bh,bhn2,h,hr1,n,nx,w}. In this case, the set of possible speeds is a half-line $[w^*,+\infty)$. Estimates of the minimal speeds $w^*$ have been derived in \cite{bhn1,bhn2,bhr,e1,e2,he,rz,z2}. Since the seminal paper of Aronson and Weinberger \cite{aw} in the homogeneous setting in $\R^N$, much work has also been devoted to asymptotic spreading speeds in KPP-type reaction-diffusion equations with compactly supported initial conditions in periodic or more general media \cite{bhn,bhn3,lyz,lz,w}, with exponentially decaying initial conditions \cite{bmr,b,hn,hr1,k,l,mr,r3,u} or with slowly decaying initial conditions \cite{cr,hr2}.

Let us now come back to the Cauchy problem (\ref{problem}) under assumption (\ref{comb}). As already emphasized, the main goal of this paper is to consider (\ref{problem}) with a very large class of front-like initial conditions, satisfying (\ref{lim12}), which are not required to converge to any constant in the direction of propagation or to be close to any pulsating front. Before stating our main results, we just need to introduce a few more notations. We consider the following linear advection-diffusion equation with the same initial data $u_0$ as for the nonlinear Cauchy problem (\ref{problem}), but with zero right-hand side:
\begin{equation}\label{heat}
\left \{
\begin{array}[c]{rcll}
v_t-\nabla\cdot(A(z)\nabla v)+q(z)\cdot\nabla v & = & 0, & t>0,\ \mbox{$z\in\overline{\Omega}$},\vspace{5pt}\\
\nu A\nabla v & = & 0, & t>0,\ \mbox{$z\in\partial\Omega$},\vspace{5pt}\\
v(0,z) & = & u_0(z), & \mbox{$z\in \Omega$},
\end{array} \right . 
\end{equation}  
Then, we introduce the following quantities, which will play an important role in the sequel 
$$\alpha_{\min}(u_0)=\lim_{t \rightarrow +\infty} \Big(\liminf_{x  \rightarrow +\infty} v(t,z) \Big)$$
and 
$$\alpha_{\max}(u_0)=\lim_{t \rightarrow +\infty} \Big(\limsup_{x  \rightarrow +\infty} v(t,z) \Big).$$
The limits in time in the previous two quantities are well-defined real numbers since the maps $t\mapsto\liminf_{x  \rightarrow +\infty}\,v(t,z)$ and $t\mapsto\limsup_{x  \rightarrow +\infty}\,v(t,z)$ are bounded (in $[0,1]$) and respectively nondecreasing and nonincreasing in time by the parabolic maximum principle (see Remark~\ref{remHeat} after Lemma~\ref{tempHeat} for more details). Furthermore, there holds
$$0 \leq \alpha_{\min}(u_0) \leq \alpha_{\max}(u_0) \leq \limsup_{x  \rightarrow +\infty} u_0(z) <\theta.$$
In particular, if $u_0(z)\to\gamma$ as $x\to+\infty$ for some real number $\gamma\in[0,\theta)$, then $\alpha_{\min}(u_0)=\alpha_{\max}(u_0)=\gamma$.

We can now state the main results of this paper. The first theorem provides lower and upper bounds for the lower and upper spreading speeds $c_*(u_0)$ and $c^*(u_0)$.  

\begin{theo}\label{mainTh}
Let $u$ be the solution of \eqref{problem} with any uniformly continuous initial condition $u_0:\Omega\to[0,1]$ satisfying \eqref{lim12}. Then
$$\displaystyle{\liminf_{t \rightarrow +\infty}}\Big(\inf_{z\in \overline{\Omega}} u(t,z)\Big)\geq \alpha_{\min} (u_0)$$
and
\be\label{ineqsc*}
c_{\alpha_{\min}(u_0)} \leq c_*(u_0)\leq c^*(u_0) \leq c_{\alpha_{\max}(u_0)},
\ee
where we recall that $c_{\alpha_{\min}(u_0)}$ and $c_{\alpha_{\max}(u_0)}$ denote the unique speeds of the pulsating fronts of $(\ref{per})$ connecting $\alpha_{\min}(u_0)$ and $\alpha_{\max}(u_0)$ to $1$. Furthermore, for every $c > c_{\alpha_{\max}(u_0)}$, there holds
$$\limsup_{ t \rightarrow +\infty}\Big(\sup_{z\in\overline{\Omega},\ \!x \geq ct} u(t,z)\Big)\leq \alpha_{\max}(u_0). $$
\end{theo}

Thus, Theorem \ref{mainTh} provides bounds for the asymptotic spreading speeds~$c_*(u_0)$ and~$c^*(u_0)$. The following theorem states a complete characterization of these sprea\-ding speeds  when the initial condition is assumed to be asymptotically periodic in the right direction.    

\begin{theo}\label{periodicTh}
Let $u_0:\Omega\to[0,1]$ be a uniformly continuous function such that there exists a uniformly continuous periodic function $w_0:\Omega\to[0,\theta)$ satisfying 
$$\lim_{ x \rightarrow +\infty }|u_0(z)-w_0(z)|=0. $$  
Then $\alpha_{\min}(u_0)=\alpha_{\max}(u_0)=\ \!<\!w_0\!>$, and consequently
$$c_*(u_0)=c^*(u_0)=c_{<w_0>},$$
where
$$<\!w_0\!>=\int_{C}\!\!\!\!\!\!\!-\ \!w_0\in[0,\theta)$$
denotes the average of the periodic function $w_0$ and $c_{<w_0>}$ is the unique speed of the pulsating travelling front of \eqref{per} connecting $<\!w_0\!>$ to $1$.
\end{theo} 

As a consequence of Theorem~\ref{periodicTh}, we have convergence in speed in the following sense: for any given value $\lambda\in(<\!w_0\!>,1)$, the set
$$E_{t,\lambda}=\{(x,y)\in\overline{\Omega},\ u(t,x,y)=\lambda\}$$
is not empty for large $t$ and
$$\lim_{t\to+\infty}\frac{\min\{x\ |\ \exists\,y,\, (x,y)\in E_{t,\lambda}\}}{t}=\lim_{t\to+\infty}\frac{\max\{x\ |\ \exists\,y,\, (x,y)\in E_{t,\lambda}\}}{t}=c_{<w_0>}.$$
In particular, for any family $(x(t),y(t))$ in $\overline{\Omega}$ such that $u(t,x(t),y(t))=\lambda$, then $x(t)/t$ converges to $c_{<w_0>}$ as $t\to+\infty$. This corresponds exactly to the notion of convergence in speed. However, it does not mean that $x(t)-c_{<w_0>}t$ converges as $t\to+\infty$ or is even bounded. But we conjecture that $x(t)-c_{<w_0>}t$ converges as $t\to+\infty$ provided that $u_0$ converges to $w_0$ sufficiently fast (exponentially) as $x\to+\infty$. A remaining open question is the convergence in profile of $u(t,\cdot)$ to the one-parameter family of time shifts of the pulsating front $u_{<w_0>}$.

In the previous two theorems, we established some general properties and bounds of the lower and upper spreading speeds $c_*(u_0)$ and $c^*(u_0)$, and we considered an important class of initial conditions for which these two quantites are equal. In what follows, we exhibit a class of initial conditions $u_0$ for which $c_*(u_0)<c^*(u_0)$ and, among other things, we will see that the behaviour of the solution $u$ along the rays with speeds $c$ between $c_*(u_0)$ and~$c^*(u_0)$ is rather complex. For the sake of clarity of the presentation, we only consider here a simple one-dimensional and homogeneous framework --more general heterogeneous equations with the same type of initial conditions and the same type of long-time behaviour could be dealt with. Consider the following Cauchy problem
\be\label{cauchy1d}\left\{\baa{rcll}
u_t-u_{xx} & = & f(u), & t>0,\ \ x\in\R,\vspace{5pt}\\
u(0,x) & = & u_0(x), & x\in\R,\eaa\right.
\ee
where the nonlinearity $f$ is of the combustion type (\ref{comb}), as above, but it depends on $u$ only. For each $\gamma\in(-\infty,\theta)$, let~$\varphi_{\gamma}$ denote the unique (up to shifts) travelling front of~(\ref{cauchy1d}) connecting $\gamma$ to $1$, with unique speed $c_{\gamma}$, that is
\be\label{cgamma}\left\{\baa{l}
\varphi_{\gamma}''+c_{\gamma}\varphi_{\gamma}'+f(\varphi_{\gamma})=0\ \hbox{ in }\R,\vspace{5pt}\\
\varphi_{\gamma}(-\infty)=1>\varphi_{\gamma}(x)>\varphi_{\gamma}(+\infty)=\gamma\ \hbox{ for all }x\in\R.\eaa\right.
\ee

\begin{theo}\label{th1d} Let $\alpha<\beta$ be any given real numbers in $[0,\theta)$. There are initial conditions $u_0:\R\to[\alpha,1]$ such that
$$\left\{\baa{l}
\displaystyle{\mathop{\liminf}_{x\to-\infty}}\,u_0(x)>\theta,\vspace{5pt}\\
\exists\,A\in\R,\ \forall\,x\ge A,\quad\alpha\le u_0(x)\le\beta\eaa\right.$$
and such that, under the general previous notations,
$$\alpha_{\min}(u_0)=\alpha,\quad\alpha_{\max}(u_0)=\beta,\quad c_*(u_0)=c_{\alpha}<c_{\beta}=c^*(u_0).$$
Furthermore,
$$\forall\,t\ge 0,\quad\alpha\le u(t,\cdot)\le 1,\quad u(t,-\infty)=1,\quad \liminf_{x\to+\infty}u(t,x)=\alpha<\beta=\limsup_{x\to+\infty}u(t,x)$$
and
$$\left\{\baa{ll}
\forall\ c<c_*(u_0),\ \forall\ A\in\R, & u(t,ct+\cdot)\displaystyle{\mathop{\longrightarrow}_{t\to+\infty}}1\ \hbox{ uniformly in }(-\infty,A],\vspace{5pt}\\
\forall\ x\in\R,\ \exists\ \alpha_x\in(\alpha,1], & \left\{\displaystyle{\mathop{\lim}_{t_k\to+\infty}}u(t_k,c_*(u_0)t_k+x)\right\}=[\alpha_x,1],\vspace{5pt}\\
\forall\ c\in(c_*(u_0),c^*(u_0)),\ \forall\ x\in\R, & \left\{\displaystyle{\mathop{\lim}_{t_k\to+\infty}}u(t_k,ct_k+x)\right\}=[\alpha,1],\vspace{5pt}\\
\forall\ x\in\R,\ \exists\ \beta_x\in[\beta,1), & \left\{\displaystyle{\mathop{\lim}_{t_k\to+\infty}}u(t_k,c^*(u_0)t_k+x)\right\}=[\alpha,\beta_x],\vspace{5pt}\\
\forall\ c>c^*(u_0),\ \forall\ A\in\R, & \displaystyle{\mathop{\lim}_{t\to+\infty}}\Big(\displaystyle{\mathop{\sup}_{x\in[A,+\infty)}}u(t,ct+x)\Big)=\beta,\vspace{5pt}\\
\forall\ c>c^*(u_0),\ \forall\ x\in\R, & \left\{\displaystyle{\mathop{\lim}_{t_k\to+\infty}}u(t_k,ct_k+x)\right\}=[\alpha,\beta].\eaa\right.$$
\end{theo}

Let us now comment the construction and the long-time behaviour of the solutions~$u$ given in Theorem~\ref{th1d}. The initial conditions~$u_0$ are constructed so that~$u_0(x)$ oscillates between~$\alpha$ and~$\beta$ as $x\to+\infty$, on larger and larger intervals. This way, the solution $u$ will somehow oscillate at large times between two approximated fronts whose speeds are approximately equal to~$c_*(u_0)=c_{\alpha}$ and~$c^*(u_0)=c_{\beta}$. In other words, the ``location"~$\xi(t)$ of the solution, that is $\xi(t)\in\R$ such that $u(t,\xi(t))=\theta$, oscillates between~$c_{\alpha}t$ and~$c_{\beta}t$, which means in particular no convergence in speed. We nevertheless provide quantitative estimates on $\xi(t)$ over some reasonably large time intervals (precise statements will be given in Section~\ref{sec1d}, see in particular the proof of Lemma~\ref{speeds} and Remark~\ref{remspeeds} below). Thus, the values of~$u(t,ct+x)$ along any ray with a given speed~$c\in(c_*(u_0),c^*(u_0))$ describe at the limit the whole interval $[\alpha,1]$, in the sense that the set of limit values of the function $t\mapsto u(t,ct+x)$ as $t\to+\infty$ is equal to the whole interval $[\alpha,1]$. In the moving frame with speed~$c_*(u_0)$ (resp.~$c^*(u_0)$), as we shall see in Section~\ref{sec1d}, the solution~$u$ is actually separated from~$\alpha$ (resp.~$1$) uniformly in~$(-\infty,A]$ (resp. $[A,+\infty)$) for any $A\in\R$, in the sense that
\be\label{c**}\left\{\baa{l}
\alpha<\displaystyle{\mathop{\liminf}_{t\to+\infty}}\Big(\displaystyle{\mathop{\inf}_{x\in(-\infty,A]}}u(t,c_*(u_0)t+x)\Big),\vspace{5pt}\\
\displaystyle{\mathop{\limsup}_{t\to+\infty}}\Big(\displaystyle{\mathop{\sup}_{x\in[A,+\infty)}}u(t,c^*(u_0)t+x)\Big)<1.\eaa\right.
\ee
However, these limits are never uniform in space, since $\inf_{\R}u(t,\cdot)=\alpha$ and $\sup_{\R}u(t,\cdot)=1$ for all $t\ge 0$.

Under the general notations and assumptions of this paper, the speed $c_*(u_0)$ is by definition the largest speed for which the solution~$u$ converges to~$1$ in any right-moving frame with a speed {\it smaller} than~$c_*(u_0)$ (and even uniformly in any given set $\{x\le A\}$). However, one of the main interests of Theorem~\ref{th1d} is to show that, even for the homogeneous equation~(\ref{cauchy1d}), the solution~$u$ may not in general be separated from~$1$ in all moving frames with speeds {\it larger} than $c_*(u_0)$: indeed, in Theorem~\ref{th1d}, there holds $\limsup_{t\to+\infty}u(t,ct+x)=1$ for all $c\in[c_*(u_0),c^*(u_0))$ and $x\in\R$, with $c_*(u_0)<c^*(u_0)$. On the other hand, again by virtue of our general definitions, the solution~$u$ is always separated from~$1$ in any right-moving frame with a speed {\it larger} than $c^*(u_0)$. However, in any such moving frame, the solution~$u$ may still have a complex behaviour and it may not converge locally to a constant in general, as seen as a byproduct of the last assertion in Theorem~\ref{th1d}. Lastly, the solution $u$ of (\ref{problem}) may not in general be separated from its infimum value in all moving frames with speeds {\it smaller} than $c^*(u_0)$, since, in the example given in Theorem~\ref{th1d}, there holds $\liminf_{t\to+\infty}u(t,ct+x)=\alpha$ for all $c\in(c_*(u_0),c^*(u_0)]$ and $x\in\R$, with $c_*(u_0)<c^*(u_0)$.

Lastly, we mention that non-convergence results similar to the ones described in Theorem~\ref{th1d} for $c\in(c_*(u_0),c^*(u_0))$ are also known to hold for the heat equation, see \cite{ce} and Remark~\ref{remNL} below. Other complex behavior may also occur for the nonlinear equation $u_t=\Delta u+f(u)$ in $\R^2$ with bistable-type nonlinearity $f$ and some appropriate initial conditions which are trapped between two shifts of a given conical front (the solutions may not in general converge to a unique shift of the given front, see~\cite{rrm}), as well as for supercritical semilinear heat equations with some initial conditions which are trapped between two ordered stationary states (the solutions may not in general converge to a unique stationary state, see \cite{py}).

%%%%%%%%%%%%%%%%%%%%%%%%%%%%%%%%%%%%%%%%
%%%%%%%%%%%%%%%%%%%%%%%%%%%%%%%%%%%%%%%%

\section{General properties}

This section is concerned with the proof of the general properties of the lower and upper spreading speeds $c_*(u_0)$ and $c^*(u_0)$. We begin in Subsection~\ref{sec21} with the proof of Theorem~\ref{periodicTh}, since it follows straightforwardly from Theorem~\ref{mainTh}. Then Subsection~\ref{sec22} is devoted to the proof of Theorem~\ref{mainTh}.

%%%%%%%%%%%%%%%%%%%%%%%%%%%%%%%%%%%%%%%%

\subsection{Asymptotically periodic initial conditions}\label{sec21}

In this subsection, we prove Theorem \ref{periodicTh}, assuming the conclusion of Theorem~\ref{mainTh}. To this end, we shall use the following theorem, providing Gaussian estimates (see Theorem~$6.1$ in \cite{daners}) for the fundamental solution of the linear equation (\ref{heat}).
 
\begin{theo}\label{GB} {\rm{\cite{daners}}}
Let $p(t,z,z')$ be the kernel of the operator $\partial_t -\nabla \cdot (A\nabla )+q \cdot \nabla$ in~$\Omega$ with no-flux boundary conditions $\nu A\nabla p=0$ on $\partial\Omega$. Then there exist some constants $C_0>0$, $\omega_1 \geq 0$ and $\omega_2 >0$ such that for all $0 <t<+\infty$ and $(z,z') \in\overline{\Omega}$
\begin{equation}\label{gaussian}
|p(t,z,z')| \leq C_0\ \!t^{-\frac{N}{2}} e^{\omega_1 t -\frac{|z-z'|^2}{\omega_2 t}}
\end{equation}
\end{theo}

In order to prove Theorem \ref{periodicTh}, we will need the following lemma, which is a consequence of the above Gaussian estimates.  
 
\begin{lem}\label{tempHeat}
Let $v$ $($resp. $w)$ be the unique solution of the linear equation \eqref{heat} in $\Omega$ with a uniformly continuous and bounded initial condition $v_0: \Omega \rightarrow\R$ $($resp. $w_0: \Omega \rightarrow\R)$. Assume furthermore that
$$\lim_{ x \rightarrow +\infty}|v_0(z)-w_0(z)|=0.$$ 
Then, for all $t \geq 0$,
$$\lim_{ x \rightarrow +\infty}|v(t,z)-w(t,z)|=0.$$ 
\end{lem}

\noindent{\bf{Proof.}} The proof uses standard arguments. We just do it here for the sake of completeness. By uniqueness of the solution of the Cauchy problem, the function $\varphi(t,z)=v(t,z)-w(t,z)$ satisfies 
$$\varphi(t,z)=\int_\Omega p(t,z,z')\,\varphi(0,z')\, dz'$$
for all $t>0$ and $z\in\overline{\Omega}$. Let $t>0$ and $\epsilon>0$ be any two arbitrary positive real numbers. From the assumption of the lemma, there is $A\in\R$ such that $|\varphi(0,z)|\le\epsilon$ for all $z=(x,y)\in\Omega$ such that $x\ge A$. Set
$$\Omega^+=\{z'=(x',y')\in\Omega,\ x'\ge A\}\ \hbox{ and }\ \Omega^-=\{z'=(x',y')\in\Omega,\ x'\le A\}.$$
For all $z=(x,y)\in\overline{\Omega}$, there holds
$$\baa{rcl}
|\varphi(t,z)| & \le & \epsilon\displaystyle{\int_{\Omega^+}}p(t,z,z')\,dz'+\displaystyle{\int_{\Omega^-}}p(t,z,z')\,|\varphi(0,z')|\,dz'\vspace{5pt}\\
& \le & \epsilon+C_0\,\|v_0-w_0\|_{\infty}\,t^{-\frac{N}{2}}\,e^{\omega_1t}\displaystyle{\int_{\Omega^-}}e^{-\frac{|z-z'|^2}{\omega_2t}}\,dz'\eaa$$
from Theorem~\ref{GB}. From (\ref{omega}), it follows that
\be\label{CNR}
|\varphi(t,z)|\le\epsilon+C_0\,\|v_0-w_0\|_{\infty}\,C_{N,R}\,t^{-\frac{N}{2}}\,e^{\omega_1t}\int_{-\infty}^{A}e^{-\frac{|x-x'|^2}{\omega_2t}}dx',
\ee
where $C_{N,R}>0$ denotes the Lebesgue measure of any euclidean ball of radius $R$ in $\R^{N-1}$. Since the last integral does not depend on $y$ and converges to $0$ as $x\to+\infty$, one concludes that there exists $B\in\R$ such that $|\varphi(t,z)|\le 2\epsilon$ for all $z=(x,y)\in\overline{\Omega}$ such that $x\ge B$, which gives the desired conclusion. \hfill$\Box$

\begin{rem}\label{remHeat}{\rm Under the notations of Lemma~\ref{tempHeat}, it follows that the quantities $m(t)=\liminf_{x\to+\infty}v(t,z)$ and $M(t)=\limsup_{x\to+\infty}v(t,z)$ are respectively nondecreasing and nonincreasing with respect to time $t$. It is obviously sufficient to deal with $m(t)$. To do so, let $0\le t_1<t_2<+\infty$ be fixed, and let $\epsilon>0$ be arbitrary. There exist then $A\in\R$ and a uniformly continuous and bounded function $w_0:\Omega\to\R$ such that $v(t_1,\cdot)\ge w_0$ in~$\Omega$ and $w_0(z)=m(t_1)-\epsilon$ for all $z=(x,y)\in\Omega$ such that $x\ge A$. The maximum principle yields $v(t+t_1,\cdot)\ge w(t,\cdot)$ in $\overline{\Omega}$ for all $t>0$, where $w$ denotes the solution of (\ref{heat}) with initial condition $w_0$. But Lemma~\ref{tempHeat} implies that $\lim_{x\to+\infty}w(t,z)=m(t_1)-\epsilon$ for all $t\ge 0$. Consequently, $m(t_2)\ge m(t_1)-\epsilon$. The conclusion follows.}
\end{rem}

\noindent{\bf{Proof of Theorem \ref{periodicTh}.}} Let $w$ be the solution of \eqref{heat} with an initial datum $w_0$ as in the theorem. It is then classical to check that
\be\label{moyenne}
\lim_{t \rightarrow +\infty} w(t,z)=\ <\!w_0\!>\ \hbox{ uniformly in }\overline{\Omega},
\ee
where $<\!w_0\!>$ denotes the average of the periodic function $w_0$. Indeed, $w(t,\cdot)$ remains periodic for each $t>0$ by uniqueness of the Cauchy problem. Furthermore, $\min_{\overline{\Omega}}w(t,\cdot)=\min_{\overline{C}}w(t,\cdot)$ and $\max_{\overline{\Omega}}w(t,\cdot)=\max_{\overline{C}}w(t,\cdot)$ are bounded, and respectively nondecreasing and nonincreasing in $t>0$. Let $(t_n)_{n\in\N}$ be a sequence of positive times converging to $+\infty$ and $(z_n)_{n\in\N}$ be a sequence of points in $\overline{C}$ such that
$$\lim_{n\to+\infty}w(t_n,z_n)=\lim_{t\to+\infty}\Big(\min_{\overline{\Omega}}w(t,\cdot)\Big)=:m.$$
From standard parabolic estimates, up to extraction of subsequence, the $x$-periodic functions
$$w_n(t,z)=w(t+t_n,z)$$
converge locally in $t$ and uniformly in $\overline{\Omega}$ as $n\to+\infty$ to a classical solution $w_{\infty}$ of the same equation~\eqref{heat} in $\R\times\overline{\Omega}$. Furthermore, $w_{\infty}\ge m$ in $\R\times\overline{\Omega}$ and $\min_{\overline{\Omega}}w_{\infty}(0,\cdot)=m$. Thus, $w_{\infty}\equiv m$ in  $\R\times\overline{\Omega}$ from the strong parabolic maximum principle. This implies that, given any $\epsilon>0$, there is $N\in\N$ such that $|w(t_N,\cdot)-m|\le\epsilon$ in $\overline{\Omega}$, whence
$$|w(t,z)-m|\le\epsilon\ \hbox{ for all }(t,z)\in[t_N,+\infty)\times\overline{\Omega}$$
from the maximum principle. As a consequence, $w(t,z)\to m$ as $t\to+\infty$ uniformly in~$\overline{\Omega}$. On the other hand, integrating the equation (\ref{heat}) in $C$ at any time $t>0$ implies that the function $t\mapsto h(t)=\int_Cw(t,z)\,dz$ is constant in $t>0$, because $q$ is divergence-free in $\Omega$ and tangential on $\partial\Omega$. Since $w(t,z)\to w_0(z)$ as $t\to 0^+$ for all $z\in C$ and the function $w$ is globally bounded, Lebesgue's dominated convergence theorem implies that $h(t)=\,<\!w_0\!>|C|$ for all $t>0$, where $|C|$ denotes the Lebesgue measure of the periodicity cell $C$. Eventually, this yields~(\ref{moyenne}).\par
Therefore, by the uniformity of the limit (\ref{moyenne}) and by Lemma \ref{tempHeat}, we deduce that 
$$\alpha_{\min}(u_0)=\alpha_{\max}(u_0)=<\!w_0\!>.$$
But $c_{\alpha_{\min}(u_0)} \leq c_*(u_0) \leq c^*(u_0) \leq c_{\alpha_{\max}(u_0)}$ from Theorem \ref{mainTh}. Hence,
$$c_*(u_0)=c^*(u_0)=c_{<w_0>},$$
and the proof of Theorem~\ref{periodicTh} is complete.\hfill$\Box$

%%%%%%%%%%%%%%%%%%%%%%%%%%%%%%%%%%%%%%%%

\subsection{Lower and upper bounds for $c_*(u_0)$ and $c^*(u_0)$}\label{sec22}

The following section is devoted to the proof of Theorem~\ref{mainTh}. We start with a general result ensuring that any solution of (\ref{heat}) with a compactly supported initial datum converges uniformly to $0$ as $t\to+\infty$. 

\begin{lem}\label{mass}
Let $w$ be the solution of the linear equation
\be\label{heatw}\left\{\begin{array}[c]{rcll}
w_t-\nabla\cdot(A(z)\nabla w)+q(z)\cdot\nabla w & = & 0, & t>0,\ \mbox{$z\in\overline{\Omega}$},\vspace{5pt}\\
\nu A\nabla w & = & 0, & t>0,\ \mbox{$z\in\partial\Omega$},\vspace{5pt}\\
w(0,z) & = & w_0(z), & \mbox{$z\in \Omega$},\end{array}\right.
\ee
where $w_0:\overline{\Omega}\to\R$ is continuous and compactly supported in $\overline{\Omega}$. Then, for all $t>0$,
$$\int_{\Omega} w(t,z)\,dz=\int_{\Omega}w_0(z)\,dz.$$
Furthermore, $w(t,z)\to 0$ as $t\to+\infty$ uniformly in $\overline{\Omega}$.
\end{lem}

\noindent{\bf{Proof.}} First of all, since $w_0$ has a compact support, denoted by $K=\hbox{supp}(w_0)$, the following pointwise estimate follows from Theorem \ref{GB}:
\be\label{supw}
\forall\,t>0,\ \forall\,z\in\overline{\Omega},\quad|w(t,z)|\leq C_0\,\|w_0\|_{\infty}\,t^{-\frac{N}{2}}\,e^{\omega_1t}\int_Ke^{-\frac{|z-z'|^2}{\omega_2 t}}dz'.
\ee
By pointwise gradient bounds (see \cite{lsu}), we get that, for every $t>0$ and every $z\in\overline{\Omega}$,
$$|\nabla_{z}w(t,z)|\leq C(t)\max_{t'\in[t/2,t],\, z'\in\overline{\Omega},\,|z'-z|\le 1}|w(t',z')|,$$
where $C(t)$ depends on $t$ but not on $z$. As a consequence, for any $0<a\le b<+\infty$, there are positive constants $C'_{a,b}$ and $\omega_{a,b}$ which depend on $a$ and $b$, such that
$$\forall\,t\in[a,b],\ \forall\,z\in\overline{\Omega},\quad |w(t,z)|+|\nabla_zw(t,z)|\le C'_{a,b}\,e^{-\omega_{a,b}|z|^2}.$$
Notice in particular that the integrals of $w$ and $|\nabla_zw|$ over $\Omega$ converge at any time $t>0$.\par
Fix now any two times $0<t<t'$. Integrate the equation
$$w_t=\nabla\cdot(A(z)\nabla w)-q\cdot\nabla w$$
over $[t,t']\times(\Omega\cap B_R)$, where $B_R$ denotes the euclidean ball of $\R^N$ centered at the origin with radius $R$, and pass to the limit as $R\to+\infty$. It follows from the previous estimates that
$$\int_{\Omega}w(t,z)\,dz=\int_{\Omega}w(t',z)\,dz,$$
using once again the assumptions that $q$ is divergence-free in $\Omega$ and tangential on the boundary $\partial\Omega$.\par
Lastly, we know that $w(t,z)\to w_0(z)$ as $t\to 0^+$ for all $z\in\Omega$. Moreover, $|w|\le\|w_0\|_{\infty}$ in $[0,+\infty)\times\Omega$. On the other hand, there is $\eta>0$ such that
$$|z-z'|\ge\eta\,|z|\ \hbox{ for all }z'\in K\hbox{ and }z\in\overline{\Omega}\ \backslash\ (K+B_1).$$
Therefore, it follows from (\ref{supw}) that, for all $0<t\le\min(1,1/\omega_2)$ and $z\in\overline{\Omega}\ \backslash\ (K+B_1)$,
$$|w(t,z)|\le C_0\,\|w_0\|_{\infty}\,\omega_2^{\frac{N}{2}}\,e^{\omega_1}\int_{\frac{K-z}{\sqrt{\omega_2t}}}e^{-|z''|^2}dz''\le C_0\,\|w_0\|_{\infty}\,\omega_2^{\frac{N}{2}}\,e^{\omega_1}\int_{\R^N\backslash B_{\eta|z|}}e^{-|z''|^2}dz''.$$
Since the right-hand side does not depend on $t$ and is integrable (with respect to $z$) over $\overline{\Omega}\ \backslash\ (K+B_1)$, Lebesgue's dominated convergence theorem finally yields that
$$\int_{\Omega}w(t,z)\,dz\to\int_{\Omega}w_0(z)\,dz\ \hbox{ as }t\to 0^+.$$
Hence, for every $t>0$, the integral of $w(t,\cdot)$ over $\Omega$ is the same as that of $w_0$.\par
Let us now prove that $w$ converges to $0$ as $t\to+\infty$ uniformly in $\overline{\Omega}$. Consider first the case when $w_0$ is nonnegative, whence $w(t,\cdot)\ge 0$ in $\overline{\Omega}$ for all $t>0$. The quantity 
$$\ell(t)=\sup_\Omega w(t,\cdot)$$
belongs to $[0,\|w_0\|_{\infty}]$ and it is nonincreasing with respect to $t\ge 0$, from the maximum principle. As a consequence, it has a limit in $[0,\|w_0\|_{\infty}]$ when $t \rightarrow +\infty$, denoted $\ell_\infty$. Assume that $\ell_\infty>0$. Then there exist a sequence $(t_n)_{n\in\N}$ and a sequence of points $(z_n)_{n\in\N}$ in $\overline{\Omega}$ such that $t_n\to+\infty$ and $w(t_n,z_n)\to\ell_{\infty}$ as $n\to+\infty$. Denote $z_n=(x_n,y_n)=(k_nL+x'_n,y_n)$, where $k_n\in\Z$ and $(x'_n,y_n)\in\overline{C}$. Up to extraction of a subsequence, the points $(x'_n,y_n)$ converge to $z_{\infty}\in\overline{C}$ and the functions
$$w_n(t,z)=w_n(t,x,y)=w(t+t_n,x+k_nL,y)$$
converge locally uniformly in $\R\times\overline{\Omega}$ to a classical bounded solution $w_{\infty}$ of the same equation as $w$, such that $w_{\infty}\le\ell_{\infty}$ in $\R\times\overline{\Omega}$ and $w_{\infty}(0,z_{\infty})=\ell_{\infty}$. Therefore, $w_{\infty}\equiv\ell_{\infty}$ in $\R\times\overline{\Omega}$ from the strong parabolic maximum principle. In other words, the functions $w_n(t,z)=w_n(t,x,y)=w(t+t_n,x+k_nL,y)$ converge locally uniformly in $\R\times\overline{\Omega}$ to the positive constant~$\ell_{\infty}$, which implies that the integrals of the nonnegative functions $w(t_n,\cdot)$ over~$\Omega$ cannot stay bounded. This leads to a contradiction. Thus $w(t,\cdot)\to 0$ as $t\to+\infty$ uniformly in $\overline{\Omega}$.\par
In the general case when $w_0$ has no sign, one can write $w_0=w^+_0-w^-_0$, where $w^+_0(x)=\max(w_0(x),0)$ and $w^-_0(x)=\max(-w_0(x),0)$ for all $x\in\overline{\Omega}$. By uniqueness and linearity of the Cauchy problem (\ref{heatw}), it follows that $w(t,z)=w_1(t,z)-w_2(t,z)$ for all $t>0$ and $z\in\overline{\Omega}$, where $w_1$ and $w_2$ solve (\ref{heatw}) with initial conditions $w^+_0$ and $w^-_0$ respectively. But the previous paragraph implies that $w_1$ and $w_2$ converge to $0$ as $t\to+\infty$ uniformly in $\overline{\Omega}$, whence $\lim_{t\to+\infty}\|w(t,\cdot)\|_{\infty}=0$ and the proof of Lemma~\ref{mass} is now complete.\hfill$\Box$\break

The next lemma provides the proof of the first statement of Theorem \ref{mainTh}. 

\begin{lem}\label{L0}
Let $u$ be a solution of $\eqref{problem}$ satisfying the assumptions of Theorem $\ref{mainTh}$. Then~$u$ satisfies 
\begin{equation*}
\liminf_{ t \rightarrow +\infty}\Big(\inf_{z \in \overline \Omega} u(t,z)\Big) \geq \alpha_{\min}(u_0).
\end{equation*}
\end{lem}

\noindent{\bf{Proof.}} Fix any arbitrary $\varepsilon >0$. By definition of $\alpha_{\min}(u_0)$, there exist $T>0$ and $A_1>0$ such that for every $z=(x,y)\in\overline{\Omega}$ with $x \geq A_1$, there holds
$$v(T,z) \geq \alpha_{\min}(u_0)-\varepsilon. $$
Since $\liminf_{x\to-\infty}u_0(z)>\theta>\alpha_{\min}(u_0)$ and since the map $t\mapsto\liminf_{x\to-\infty}v(t,z)$ is nondecreasing (with the same kind of arguments as in Lemma~\ref{tempHeat} and Remark~\ref{remHeat}), it follows that there exists $A_2<0$ such that for every $z=(x,y)\in\overline{\Omega}$ with $x \leq A_2$, there holds
$$v(T,z) \geq\theta>\alpha_{\min}(u_0)-\varepsilon. $$
Denote $K$ the compact set
$$K=\{z=(x,y)\in\overline{\Omega},\ A_1\le x\le A_2\}.$$
Consequently,
$$\forall\,z\in\overline{\Omega}\backslash K,\quad v(T,z)\ge\alpha_{\min}(u_0)-\epsilon.$$
In particular, since $v$ is also globally bounded, there exists a continuous and compactly supported function $w_0:\overline{\Omega}\to\R$ such that
$$\forall\, z\in\overline{\Omega},\quad v(T,z)\ge\alpha_{\min}(u_0)-\epsilon-w_0(z).$$
By linearity and uniqueness of the Cauchy problem for (\ref{heatw}), it follows that
$$\forall\,t>0,\ \forall\,z\in\overline{\Omega},\quad v(T+t,z)\ge\alpha_{\min}(u_0)-\epsilon-w(t,z),$$
where $w$ is the solution of (\ref{heatw}) with initial condition $w_0$. On the other hand, using Lemma \ref{mass}, we know that $\|w(t,\cdot)\|_{\infty}\to 0$ as $t\to+\infty$. Hence,
$$\inf_{\overline{\Omega}}v(t,\cdot)\ge\alpha_{\min}(u_0)-2\epsilon$$
for $t$ large enough. This gives directly the desired result, letting $\varepsilon $ going to zero and using the fact that $u \geq v$ for every $t>0$ and $z \in\overline{\Omega}$ by the maximum principle (because $f$ is nonnegative).\hfill$\Box$\break

We now come to the proof of the inequalities (\ref{ineqsc*}) of Theorem~\ref{mainTh}. We first prove the following lemma, which gives directly the first inequality, namely $c_{\alpha_{\min}(u_0)} \leq c_*(u_0)$.

\begin{lem}\label{L1}
Let $u$ be a solution of $\eqref{problem}$ satisfying the assumptions of Theorem $\ref{mainTh}$. Then
\begin{equation*}
\forall\,c<c_{\alpha_{\min}(u_0)},\ \forall\,A\in\R,\quad\lim_{t \rightarrow +\infty}\Big(\inf_{x\le A,\,(x+ct,y)\in\overline{\Omega}}u(t,x+ct ,y)\Big)=1.
\end{equation*}
\end{lem}

\noindent{\bf{Proof.}} First, observe that there exist $\beta>0$, $A<0$ and a uniformly continuous function $U_0:\Omega\to[0,1]$ such that
$$U_0\le u_0\hbox{ in }\Omega\ \hbox{ and }\ U_0(x,y)=\theta+\beta\hbox{ if }x\le A.$$
Let $\underline{f}:\R\to\R_+$ be the function defined by $\underline{f}(s)=\min_{z\in\overline{\Omega}}f(z,s)$ for all $s\in\R$. Let $U$ be the solution of the Cauchy problem (\ref{problem}) with the function $f$ being replaced by $\underline{f}$, and with initial condition $U_0$. The maximum principle yields $0\le U(t,z)\le u(t,z)\le 1$ for all $t>0$ and $z\in\overline{\Omega}$. Set $\xi(t)=\lim_{x\to-\infty}U(t,z)$ for all $t\ge 0$. The function $\xi:\R_+\to[0,1]$ satisfies $\xi(0)=\theta+\beta\in(\theta,1]$ and $\xi'(t)=\underline{f}(\xi(t))$ for all $t\ge0$. From the assumption on $f$ made in (\ref{comb}), one concludes that $\xi(t)\to 1$ as $t\to+\infty$. Therefore,
\be\label{liminfu}
\liminf_{x\to-\infty}u(t,z)\to 1\ \hbox{ as }t\to+\infty.
\ee\par
Then, let us consider a family of continuous functions $(\underline{f}_{\eta})_{0\le\eta<1-\theta}:\overline{\Omega}\times\R\to\R$, with $\underline{f}_0=f$, such that each function $\underline{f}_{\eta}$ is periodic with respect to $z$, is of class $C^{0,\alpha}$ with respect to $z$ locally uniformly in $u$, has a restriction to $\overline{\Omega}\times[0,1-\eta]$ which is of class $C^1$ with respect to $u$, and satisfies
$$\forall\,z\in\overline{\Omega},\left\{\begin{array}{l}
\underline{f}_{\eta}(z,\cdot) \equiv 0 \,\,\,\mbox{on}\,\, (-\infty,\theta] \cup
[1-\eta,+\infty),\vspace{5pt}\\
\underline{f}_{\eta}(z,\cdot)>0\,\,\,\,\mbox{on}\,\,(\theta,1-\eta),\,\,\,\,\displaystyle{\frac{\partial\underline{f}_{\eta}}{\partial u}}(z,(1-\eta)^-)=-\displaystyle{\mathop{\lim}_{s\to 0^+}}\displaystyle{\frac{\underline{f}_{\eta}(z,1-\eta-s)}{s}}<0.\end{array}\right.$$
Furthermore, the functions $(\underline{f}_{\eta})$ are chosen in such a way that $\eta\mapsto\underline{f}_{\eta}(z,u)$ is nonincreasing in $[0,1-\theta)$ for each $(z,u)\in\overline{\Omega}\times\R$, and $\underline{f}_{\eta}\to f$ as $\eta\to0$ uniformly in $\overline{\Omega}\times\R$. For each $\gamma\in(-\infty,\theta)$ and $\eta\in[0,1-\theta)$, it is known \cite{bh} that there exists a unique speed $\underline{c}_{\gamma,\eta}>0$ and a unique (up to time shifts) front
$$\R\times\overline{\Omega}\ni(t,z)\mapsto\underline{u}_{\gamma,\eta}(t,z)=\underline{\varphi}_{\gamma,\eta}(x-\underline{c}_{\gamma,\eta}t,z)\ \in(\gamma,1-\eta)$$
solving (\ref{per}) with $\underline{f}_{\eta}$ instead of $f$, and satisfying (\ref{ptf}) and (\ref{limits}) with $\underline{c}_{\gamma,\eta}$ and $1-\eta$ instead of $c$ and $1$, respectively. From Theorem~\ref{thFront}, the speeds $c_{\gamma}=\underline{c}_{\gamma,0}$ are continuous with respect to $\gamma<\theta$. It also follows from \cite{bh,bhbook} that $\underline{c}_{\gamma,\eta}\to\underline{c}_{\gamma,0}=c_{\gamma}$ as $\eta\to 0$, for each $\gamma\in(-\infty,\theta)$.\par
Fix now any real number $c$ such that $c<c_{\alpha_{\min}(u_0)}$, any real number $A$ and any positive real number $\epsilon>0$. From the previous paragraph, one can choose $\kappa>0$ small enough and then $\eta\in(0,\epsilon)$ small enough so that
$$c<\underline{c}_{\alpha_{\min}(u_0)-\kappa,\eta}=:c'.$$\par
In order to conclude, we will put below the solution $u$ of (\ref{problem}) a pulsating front subsolution which will travel at speed $\underline{c}_{\alpha_{\min}(u_0)-\kappa,\eta}$ and will be larger than $1-\epsilon$ on the left. Indeed, from Lemma~\ref{L0} and (\ref{liminfu}), there exists a time $T>0$ such that
$$\liminf_{x\to-\infty}u(T,z)>1-\eta\ \hbox{ and }\ \inf_{z\in\overline{\Omega}}u(T,z)>\alpha_{\min}(u_0)-\kappa.$$
Since $\underline{\varphi}_{\alpha_{\min}(u_0)-\kappa,\eta}(-\infty,\cdot)=1-\eta$ and $\underline{\varphi}_{\alpha_{\min}(u_0)-\kappa,\eta}(+\infty,\cdot)=\alpha_{\min}(u_0)-\kappa$, there exists then a time-shift $T_0\in\R$ such that
$$u(T,z)\ge\underline{u}_{\alpha_{\min}(u_0)-\kappa,\eta}(T+T_0,z)=\underline{\varphi}_{\alpha_{\min}(u_0)-\kappa,\eta}(x-c'(T+T_0),z)\ \hbox{ for all }z\in\overline{\Omega}.$$
Since $\underline{f}_{\eta}\le f$, the function $\underline{u}_{\alpha_{\min}(u_0)-\kappa,\eta}$ is a subsolution of the equation (\ref{per}), whence
$$u(t,z)\ge\underline{u}_{\alpha_{\min}(u_0)-\kappa,\eta}(t+T_0,z)=\underline{\varphi}_{\alpha_{\min}(u_0)-\kappa,\eta}(x-c'(t+T_0),z)\ \hbox{ for all }t\ge T\hbox{ and }z\in\overline{\Omega}$$
from the maximum principle. In particular, for all $t\ge T$,
$$\inf_{x\le A,\,(x+ct,y)\in\overline{\Omega}}u(t,x+ct,y)\ge\inf_{x\le A,\,(x+ct,y)\in\overline{\Omega}}\underline{\varphi}_{\alpha_{\min}(u_0)-\kappa,\eta}(x+(c-c')t-c'T_0,x+ct,y).$$
But since $c<c'$ and $\underline{\varphi}_{\alpha_{\min}(u_0)-\kappa,\eta}(-\infty,\cdot)=1-\eta>1-\epsilon$ uniformly in $\overline{\Omega}$, one concludes that
$$\inf_{x\le A,\,(x+ct,y)\in\overline{\Omega}}u(t,x+ct,y)\ge 1-\epsilon$$
for $t$ large enough. That completes the proof of Lemma~\ref{L1}.\hfill$\Box$\break

The next lemma is a key step which will lead to the end of the proof of Theorem~\ref{mainTh}.

\begin{lem}\label{L2}
Let $u$ be a solution of \eqref{problem} and $v$ a solution of
\eqref{heat}, with the same initial condition $u_0$ satisfying the assumptions of Theorem~$\ref{mainTh}$. Then, for all $t\geq 0$,
\begin{equation*}
\lim_{x \rightarrow +\infty}\Big(u(t,z)-v(t,z)\Big)=0. 
\end{equation*}  
\end{lem}

We postpone the proof of this technical lemma to the end of this subsection and we finish the proof of Theorem \ref{mainTh}.\hfill\break

\noindent{\bf{End of the proof of Theorem \ref{mainTh}.}} We shall prove 
\begin{equation}\label{maxu}
\limsup_{ t \rightarrow +\infty}\Big(\sup_{x\ge A,\,(x+ct,y)\in\overline{\Omega}}u(t,x+ct,y)\Big)\leq \alpha_{\max}(u_0)
\end{equation}
for any $A\in\R$ and for any speed $c$ such that $c > c_{\alpha_{\max}(u_0)}$. This will give the last assertion of Theorem~\ref{mainTh} and will also imply that $c^*(u_0)\le c_{\alpha_{\max}(u_0)}$.\par
Observe first that
\begin{equation}\label{compu}
\lim_{ t \rightarrow +\infty}\Big(\limsup_{x  \rightarrow +\infty} u(t,z)\Big)\leq \alpha_{\max}(u_0). 
\end{equation}
Indeed, the same property holds for $v$ by definition of $\alpha_{\max}(u_0)$, where $v$ is the solution of~(\ref{heat}) with the same initial condition $u_0$ as $u$. Therefore, (\ref{compu}) follows from Lemma \ref{L2}.\par
We shall then construct a pulsating travelling front which will be a supersolution for $u$ and which will travel to the right at a speed larger than but close to $c_{\alpha_{\max}(u_0)}$. The proof proceeds in a similar way as in Lemma~\ref{L1}. Consider a family of continuous functions $(\overline{f}_{\eta})_{\eta\ge 0}:\overline{\Omega}\times\R\to\R$, with $\overline{f}_0=f$, such that each function $\overline{f}_{\eta}$ is periodic with respect to~$z$, is of class $C^{0,\alpha}$ with respect to~$z$ locally uniformly in~$u$, has a restriction to $\overline{\Omega}\times[0,1+\eta]$ which is of class $C^1$ with respect to~$u$, and satisfies
$$\left\{\begin{array}{l}
\overline{f}_{\eta}(z,\cdot) \equiv 0 \,\,\,\mbox{on}\,\, (-\infty,\theta] \cup
[1+\eta,+\infty),\vspace{5pt}\\
\overline{f}_{\eta}(z,\cdot)>0\,\,\,\,\mbox{on}\,\,(\theta,1+\eta),\,\,\,\,\displaystyle{\frac{\partial\overline{f}_{\eta}}{\partial u}}(z,(1+\eta)^-)=-\displaystyle{\mathop{\lim}_{s\to 0^+}}\displaystyle{\frac{\overline{f}_{\eta}(z,1+\eta-s)}{s}}<0.\end{array}\right.$$
Furthermore, the functions $(\overline{f}_{\eta})$ are chosen in such a way that $\eta\mapsto\overline{f}_{\eta}(z,u)$ is nondecreasing in $[0,+\infty)$ for each $(z,u)\in\overline{\Omega}\times\R$, and $\overline{f}_{\eta}\to f$ as $\eta\to0$ uniformly in $\overline{\Omega}\times[0,1]$. For each $\gamma\in(-\infty,\theta)$ and $\eta\in[0,+\infty)$, there exists a unique speed $\overline{c}_{\gamma,\eta}>0$ and a unique (up to time shifts) front
$$\R\times\overline{\Omega}\ni(t,z)\mapsto\overline{u}_{\gamma,\eta}(t,z)=\overline{\varphi}_{\gamma,\eta}(x-\overline{c}_{\gamma,\eta}t,z)\ \in(\gamma,1+\eta)$$
solving (\ref{per}) with $\overline{f}_{\eta}$ instead of $f$, and satisfying (\ref{ptf}) and (\ref{limits}) with $\overline{c}_{\gamma,\eta}$ and $1+\eta$ instead of $c$ and $1$, respectively. Furthermore, $\gamma\mapsto\overline{c}_{\gamma,0}=c_{\gamma}$ is continuous and $\overline{c}_{\gamma,\eta}\to\overline{c}_{\gamma,0}=c_{\gamma}$ as $\eta\to 0$, for each $\gamma\in(-\infty,\theta)$. Fix now any real number $c$ such that $c>c_{\alpha_{\max}(u_0)}$, any real number $A$ and any positive real number $\epsilon>0$. One can then choose $\kappa\in(0,\epsilon)$ small enough and then $\eta>0$ small enough so that
$$c>\overline{c}_{\alpha_{\max}(u_0)+\kappa,\eta}=:c'.$$
From (\ref{compu}), there exists a time $T>0$ such that
$$\limsup_{x\to+\infty}u(T,z)<\alpha_{\max}(u_0)+\kappa.$$
Since $u$ is also such that $u(T,\cdot)\le 1$ in $\overline{\Omega}$, there exists then a time-shift $T_0\in\R$ such that
$$u(T,z)\le\overline{u}_{\alpha_{\max}(u_0)+\kappa,\eta}(T+T_0,z)=\overline{\varphi}_{\alpha_{\max}(u_0)+\kappa,\eta}(x-c'(T+T_0),z)\ \hbox{ for all }z\in\overline{\Omega}.$$
Since $\overline{f}_{\eta}\ge f$, the function $\overline{u}_{\alpha_{\max}(u_0)+\kappa,\eta}$ is a supersolution of the equation (\ref{per}), whence
$$u(t,z)\le\overline{u}_{\alpha_{\max}(u_0)+\kappa,\eta}(t+T_0,z)=\overline{\varphi}_{\alpha_{\max}(u_0)+\kappa,\eta}(x-c'(t+T_0),z)\ \hbox{ for all }t\ge T\hbox{ and }z\in\overline{\Omega}$$
from the maximum principle. In particular, for all $t\ge T$,
$$\sup_{x\ge A,\,(x+ct,y)\in\overline{\Omega}}u(t,x+ct,y)\le\sup_{x\ge A,\,(x+ct,y)\in\overline{\Omega}}\overline{\varphi}_{\alpha_{\max}(u_0)+\kappa,\eta}(x+(c-c')t-c'T_0,x+ct,y).$$
But since $c>c'$ and $\overline{\varphi}_{\alpha_{\max}(u_0)+\kappa,\eta}(+\infty,\cdot)=\alpha_{\max}(u_0)+\kappa<\alpha_{\max}(u_0)+\epsilon$ uniformly in $\overline{\Omega}$, one concludes that
$$\sup_{x\ge A,\,(x+ct,y)\in\overline{\Omega}}u(t,x+ct,y)\le\alpha_{\max}(u_0)+\epsilon$$
for $t$ large enough. That completes the proof of (\ref{maxu}) and Theorem~\ref{mainTh}.\hfill$\Box$\break

\noindent{\bf{Proof of Lemma \ref{L2}.}} From the maximum principle, we have
\begin{equation*}
0 \leq u(t,z)-v(t,z)\ \hbox{ for all }t>0\hbox{ and }z\in\overline{\Omega} 
\end{equation*} 
because $f\ge 0$. We denote $w(t,z)=u(t,z)-v(t,z)$. Notice that $w(0,z)=0$ and the conclusion of Lemma~\ref{L2} holds immediately at time $t=0$. We are now going to construct a suitable supersolution for $w$. Let $\gamma<\theta$ be such that
$$\limsup_{x\to+\infty}u_0(z)<\gamma.$$
Pick any positive real number $\eta>0$ and consider the front $\overline{u}_{\gamma,\eta}$ connecting $1+\eta$ to $\gamma$, for the nonlinearity $\overline{f}_{\eta}$, under the above notations in the proof of Theorem~\ref{mainTh}. Denote $c=\overline{c}_{\gamma,\eta}>0$ its speed. Since $u_0$ is also not larger than $1$, there exists $T_0\in\R$ such that $u_0\le\overline{u}_{\gamma,\eta}(T_0,\cdot)$ in $\Omega$. Therefore,
$$u(t,z)\le\overline{u}_{\gamma,\eta}(t+T_0,z)=\overline{\varphi}_{\gamma,\eta}(x-c(t+T_0),z)$$
for all $t>0$ and $z\in\overline{\Omega}$ from the maximum principle. Since $\overline{\varphi}_{\gamma,\eta}(+\infty,\cdot)=\gamma<\theta$, it follows that there exists a constant $D>0$ such that for all $t >0$ and all $z=(x,y)\in\overline{\Omega}$ such that $x\geq ct+D$, there holds $u(t,z) \leq \theta$.% (see Figure $1$). 
%\begin{figure}[htbp]
%\begin{center}
%\includegraphics[scale=0.3]{super.eps}
%\end{center}
%\caption{\label{multi} The supersolution for $w$ in a one-dimensional framework.}
%\end{figure} 

We now use Duhamel's formula to express the solution $w$ of the problem
$$w_t-\nabla\cdot(A(z)\nabla w)+q(z)\cdot\nabla w=f(u).$$
Denoting $\mathcal S(t)=e^{-t{\mathcal{L}}}$ the strongly continuous semi-group generated by the operator $\mathcal L= -\nabla \cdot (A(z)\nabla )+q \cdot \nabla$ with Neumann boundary conditions $\nu A(z)\nabla=0$ on $\partial\Omega$, we have 
$$w(t,z)=\int_0^t \mathcal S(t-s)[f(u(s,\cdot))](z)\,ds.$$
for all $t>0$ and $z\in\overline{\Omega}$. Therefore, with the notations of Theorem~\ref{GB}, we get
$$\forall\,t>0,\ \forall\,z\in\overline{\Omega},\quad w(t,z)=\int_0^t\!\!\int_\Omega p(t-s,z,z')f(u(s,z'))\,dz'\,ds.$$\par
Fix now any $t>0$ and $\varepsilon >0$. Choose any $\delta$ such that 
$0<\delta<\min(t,\varepsilon/\|f\|_{\infty})$ and write 
$$w(t,z)=\underbrace{\int_0^{t-\delta}\!\!\!\!\int_\Omega p(t-s,z,z')f(u(s,z'))\,dz'\,ds}_{=:I(t,z)}+\underbrace{\int_{t-\delta}^t\!\int_\Omega p(t-s,z,z')f(u(s,z'))\,dz'\,ds}_{=:I\!I(t,z)}$$
for all $z\in\overline{\Omega}$. Notice that $0\le I\!I(t,z)\le\|f\|_{\infty}\delta\le\varepsilon$. Let us now estimate the integral~$I(t,z)$. By the Gaussian estimates in Theorem~\ref{GB}, we get that 
$$0\le I(t,z) \leq C_0\int_0^{t-\delta}\!\!\!\!\int_{\Omega}(t-s)^{-\frac{N}{2}}e^{\omega_1(t-s)-\frac{|z-z'|^2}{\omega_2(t-s)}}f(u(s,z'))\,dz'\,ds$$
for all $z\in\overline{\Omega}$. Remember that, for $z'=(x',y')\in\overline{\Omega}$, there holds $u(s,z')\le\theta$, whence $f(u(s,z'))=0$, as soon as $x'\ge cs+D$. Consequently, for all $z\in\overline{\Omega}$,
$$\baa{rcl}
0\ \le\ I(t,z) & \le & C_0\,\|f\|_{\infty}\displaystyle{\int_0^{t-\delta}}\!\!\!\displaystyle{\int_{\{z'=(x',y'),\,x'<cs+D\}}}(t-s)^{-\frac{N}{2}}e^{\omega_1(t-s)-\frac{|z-z'|^2}{\omega_2(t-s)}}\,dz'\,ds\vspace{5pt}\\
& \le & C_0\,\|f\|_{\infty}\,\delta^{-\frac{N}{2}}\,e^{\omega_1t}\,C_{N,R}\displaystyle{\int_0^{t-\delta}}\!\!\!\!\displaystyle{\int_{-\infty}^{cs+D}}e^{-\frac{|x-x'|^2}{\omega_2t}}dx'\,ds,\eaa$$
where $C_{N,R}>0$ is given as in (\ref{CNR}). But the right-hand side of the last inequality does not depend on $y$ and goes to $0$ as $x\to+\infty$. It follows that $\lim_{x\to+\infty}I(t,z)=0$, whence $0\le w(t,z)\le 2\,\varepsilon$ for all $z=(x,y)\in\overline{\Omega}$ such that $x\ge B$, for some large enough $B$. Since $\varepsilon$ is arbitrary small, this gives the desired result.\hfill$\Box$

%%%%%%%%%%%%%%%%%%%%%%%%%%%%%%%%%%%%%%%%
%%%%%%%%%%%%%%%%%%%%%%%%%%%%%%%%%%%%%%%%

\SE{Example for which $c_*(u_0)<c^*(u_0)$}\label{sec1d}

This section is devoted to the proof of Theorem~\ref{th1d}. Let $f=f(u)$ be a nonlinearity satisfying~(\ref{comb}). Let $\alpha$ and $\beta$ be given throughout the section, such that
$$0\le\alpha<\beta<\theta.$$

%%%%%%%%%%%%%%%%%%%%%%%%%%%%%%%%%%%%%%%%

\subsection{Proof of Theorem~\ref{th1d}}

In this subsection, we first define some useful notations and we derive rough estimates. Then, we state a key-lemma which enables us to complete the proof of the theorem.

\subsubsection*{Approximating fronts}

As in the proofs of Lemma~\ref{L1} and Theorem~\ref{mainTh}, we consider two families $(\underline{f}_{\eta})_{\eta\in[0,1-\theta)}$ and $(\overline{f}_{\eta})_{\eta\in[0,+\infty)}$ of $C^1([0,1-\eta])$ and $C^1([0,1+\eta])$ functions such that
$$\left\{\baa{l}
\forall\ \eta\in[0,1-\theta),\quad\underline{f}_{\eta}=0\hbox{ on }[0,\theta]\cup\{1-\eta\},\ \underline{f}_{\eta}>0\hbox{ on }(\theta,1-\eta),\ \underline{f}_{\eta}'(1-\eta)<0,\vspace{5pt}\\
\forall\ \eta\in[0,+\infty),\quad\overline{f}_{\eta}=0\hbox{ on }[0,\theta]\cup\{1+\eta\},\ \overline{f}_{\eta}>0\hbox{ on }(\theta,1+\eta),\ \overline{f}_{\eta}'(1+\eta)<0,\eaa\right.$$
and $\underline{f}_{\eta}$ (resp. $\overline{f}_{\eta}$) is extended by $0$ outside the interval $[0,1-\eta]$ (resp. $[0,1+\eta]$). Furthermore, these functions are chosen in such a way that $\underline{f}_0=\overline{f}_0=f$, that
$$\left\{\baa{ll}
\underline{f}_{\eta_1}\ge\underline{f}_{\eta_2}\hbox{ in }\R & \hbox{ if }\ 0\le\eta_1\le\eta_2<1-\theta,\vspace{5pt}\\
\overline{f}_{\eta_1}\le\overline{f}_{\eta_2}\hbox{ in }\R & \hbox{ if }\ 0\le\eta_1\le\eta_2,\eaa\right.$$
and that $\lim_{\eta\to 0}\|\underline{f}_{\eta}-f\|_{L^{\infty}(\R)}=\lim_{\eta\to 0}\|\overline{f}_{\eta}-f\|_{L^{\infty}(\R)}=0$. For each $\gamma<\theta$ and $\eta\in[0,1-\theta)$, we denote $(\underline{c}_{\gamma,\eta},\underline{\varphi}_{\gamma,\eta})$ the unique solution of
\be\label{underc}\left\{\baa{l}
\underline{\varphi}_{\gamma,\eta}''+\underline{c}_{\gamma,\eta}\underline{\varphi}_{\gamma,\eta}'+\underline{f}_{\eta}(\underline{\varphi}_{\gamma,\eta})=0\ \hbox{ in }\R,\vspace{5pt}\\
\underline{\varphi}_{\gamma,\eta}(-\infty)=1-\eta>\underline{\varphi}_{\gamma,\eta}(x)>\underline{\varphi}_{\gamma,\eta}(+\infty)=\gamma\ \hbox{ for all }x\in\R,\quad\underline{\varphi}_{\gamma,\eta}(0)=\theta.\eaa\right.
\ee
Similarly, for each $\gamma<\theta$ and $\eta\in[0,+\infty)$, we denote $(\overline{c}_{\gamma,\eta},\overline{\varphi}_{\gamma,\eta})\in\R\times C^2(\R)$ the unique solution of
\be\label{overc}\left\{\baa{l}
\overline{\varphi}_{\gamma,\eta}''+\overline{c}_{\gamma,\eta}\overline{\varphi}_{\gamma,\eta}'+\overline{f}_{\eta}(\overline{\varphi}_{\gamma,\eta})=0\ \hbox{ in }\R,\vspace{5pt}\\
\overline{\varphi}_{\gamma,\eta}(-\infty)=1+\eta>\overline{\varphi}_{\gamma,\eta}(x)>\overline{\varphi}_{\gamma,\eta}(+\infty)=\gamma\ \hbox{ for all }x\in\R,\quad\overline{\varphi}_{\gamma,\eta}(0)=\theta.\eaa\right.
\ee
Notice that with the normalization of $\underline{\varphi}_{\gamma,\eta}$ and $\overline{\varphi}_{\gamma,\eta}$ at $0$, these functions are then really unique. With these notations, for each $\gamma<\theta$, there holds $\underline{c}_{\gamma,0}=\overline{c}_{\gamma,0}=c_{\gamma}$ and the functions~$\underline{\varphi}_{\gamma,0}$ and $\overline{\varphi}_{\gamma,0}$ are equal to $\varphi_{\gamma}$ up to shifts, where $(c_{\gamma},\varphi_{\gamma})$ solves~(\ref{cgamma}). Furthermore, we recall (see~\cite{bn}) that all functions $\underline{\varphi}_{\gamma,\eta}$ and $\overline{\varphi}_{\gamma,\eta}$ are decreasing in $\R$, that the speeds $\underline{c}_{\gamma,\eta}$ and $\overline{c}_{\gamma,\eta}$ are positive, that
\be\label{monospeeds}\left\{\baa{l}
\!\!(\gamma,\eta)\!\mapsto\underline{c}_{\gamma,\eta}\hbox{ is increasing w.r.t. }\!\gamma\hbox{ and decreasing w.r.t. }\!\eta\hbox{ in }\!(-\infty,\theta)\!\times\![0,1-\theta),\vspace{5pt}\\
\!\!(\gamma,\eta)\!\mapsto\overline{c}_{\gamma,\eta}\hbox{ is increasing w.r.t. }\!\gamma\hbox{ and increasing w.r.t. }\!\eta\hbox{ in }\!(-\infty,\theta)\!\times\![0,+\infty)\eaa\right.
\ee
and that
\be\label{cgammalim}
\forall\ \gamma<\theta,\quad
\lim_{(\gamma',\eta)\to(\gamma,0)}\underline{c}_{\gamma',\eta}=\lim_{(\gamma,\eta)\to(\gamma,0)}\overline{c}_{\gamma',\eta}=c_{\gamma}>0.
\ee\par
Moreover, for each $\gamma<\theta$ and $\eta\in[0,1-\theta)$, let $\underline{u}_{\gamma,\eta}$ be the solution of the Cauchy problem
\be\label{cauchy1}\left\{\baa{l}
(\underline{u}_{\gamma,\eta})_t-(\underline{u}_{\gamma,\eta})_{xx}=f(\underline{u}_{\gamma,\eta}),\quad x\in\R,\vspace{5pt}\\
\underline{u}_{\gamma,\eta}(0,x)=\left\{\baa{ll}1-\eta, & \hbox{if }x\in(-\infty,0),\\ 1-\eta-\displaystyle{\frac{(1-\eta-\gamma)x}{2}} & \hbox{if }x\in[0,2],\\ \gamma & \hbox{if }x\in(2,+\infty).\eaa\right.\eaa\right.
\ee
For each $\gamma<\theta$, let $\overline{u}_{\gamma}$ be the solution of the Cauchy problem
\be\label{cauchy2}\left\{\baa{l}
(\overline{u}_{\gamma})_t-(\overline{u}_{\gamma})_{xx}=f(\overline{u}_{\gamma}),\quad x\in\R,\vspace{5pt}\\
\overline{u}_{\gamma}(0,x)=\left\{\baa{ll}1, & \hbox{if }x\in(-\infty,0),\\ 1-\displaystyle{\frac{(1-\gamma)x}{2}} & \hbox{if }x\in[0,2],\\ \gamma & \hbox{if }x\in(2,+\infty).\eaa\right.\eaa\right.
\ee
It is known from~\cite{k1,r3} that, for each $\eta\in[0,1-\theta)$ and each $\gamma<\theta$, there exist two real numbers $\underline{x}_{\gamma,\eta}$ and $\overline{x}_{\gamma}$ such that
\be\label{conv1}
\displaystyle{\mathop{\lim}_{t\to+\infty}}\Big(\displaystyle{\mathop{\sup}_{x\in\R}}\left|\underline{u}_{\gamma,\eta}(t,x)-\varphi_{\gamma}(x-c_{\gamma}t+\underline{x}_{\gamma,\eta})\right|\Big)=0
\ee
and
$$\displaystyle{\mathop{\lim}_{t\to+\infty}}\Big(\displaystyle{\mathop{\sup}_{x\in\R}}\left|\overline{u}_{\gamma}(t,x)-\varphi_{\gamma}(x-c_{\gamma}t+\overline{x}_{\gamma})\right|\Big)=0.$$

\subsubsection*{Definition of a class of initial conditions $u_0$}

Choose any sequence $(x_n)_{n\in\N}$ of positive real numbers such that $x_0=1$, $x_{n+1}-x_n\ge 3$ for all $n\in\N$ and
\be\label{xnn1}
\frac{x_{n+1}}{x_n}\to+\infty\ \hbox{ as }n\to+\infty.
\ee
A typical example is when $x_n=n!$ for large $n$. Let $u_0:\R\to[\alpha,1]$ be the uniformly continuous function defined by
\be\label{defu0}
u_0(x)=1\ \hbox{ if }x\in(-\infty,0),\quad u_0(x)=1-\displaystyle{\frac{(1-\alpha)x}{2}}\ \hbox{ if }x\in[0,2]=[0,x_0+1]
\ee
and, for all $n\in\N$,
\be\label{defu0bis}
u_0(x)=\left\{\baa{ll}\alpha & \hbox{ if }x\in(x_{2n}+1,x_{2n+1}-1),\vspace{5pt}\\
\alpha+\displaystyle{\frac{\beta-\alpha}{2}}\ \!(x-x_{2n+1}+1) & \hbox{ if }x\in[x_{2n+1}-1,x_{2n+1}+1],\vspace{5pt}\\
\beta & \hbox{ if }x\in(x_{2n+1}+1,x_{2n+2}-1),\vspace{5pt}\\
\beta-\displaystyle{\frac{\beta-\alpha}{2}}\ \!(x-x_{2n+2}+1) & \hbox{ if }x\in[x_{2n+2}-1,x_{2n+2}+1].\eaa\right.
\ee
Let $u$ be the solution of the Cauchy problem~(\ref{cauchy1d}) with this initial condition $u_0$. Our aim is to prove that the solution~$u$ satisfies the conclusion of Theorem~\ref{th1d}. {\it In the sequel,~$u_0$ is fixed as above and, for the sake of simplicity, we drop the dependence on $u_0$ in the quantities $\alpha_{\min}(u_0)$, $\alpha_{\max}(u_0)$, $c_*(u_0)$ and $c^*(u_0)$.}

\subsubsection*{Values of $\alpha_{\min}$ and $\alpha_{\max}$}

According to the general notations of this paper, set
$$\alpha_{\min}=\liminf_{t\to+\infty}\Big(\liminf_{x\to+\infty}v(t,x)\Big)\ \hbox{ and }\ \alpha_{\max}=\limsup_{t\to+\infty}\Big(\limsup_{x\to+\infty}v(t,x)\Big),$$
where
$$v(t,x)=\frac{1}{\sqrt{4\pi t}}\int_{-\infty}^{+\infty}e^{-\frac{|x-y|^2}{4t}}u_0(y)\,dy$$
is the solution of the heat equation
\be\label{heatv}\left\{\baa{l}
v_t=v_{xx},\quad x\in\R,\vspace{5pt}\\
v(0,\cdot)=u_0\eaa\right.
\ee
with initial condition $u_0$. Notice that $\alpha<v(t,x)<1$ for all $t>0$ and $x\in\R$, from the strong parabolic maximum principle. Observe also that, for each $t>0$,
$$v(t,x)-\beta=\frac{1}{\sqrt{4\pi t}}\int_{-\infty}^{+\infty}e^{-\frac{|x-y|^2}{4t}}(u_0(y)-\beta)\,dy\le\frac{1}{\sqrt{4\pi t}}\int_{-\infty}^2e^{-\frac{|x-y|^2}{4t}}\,dy,$$
whence $\limsup_{x\to+\infty}v(t,x)\le\beta$. For each $n\in\N$, set $y_n=\frac{x_n+x_{n+1}}{2}$. There holds
$$v(t,y_{2n})-\alpha\le\frac{1}{\sqrt{4\pi t}}\int_{|y|\ge\frac{x_{2n+1}-x_{2n}-2}{2}}e^{-\frac{|y_{2n}-y|^2}{4t}}\,dy$$
and
$$v(t,y_{2n+1})-\beta\ge-\frac{1}{\sqrt{4\pi t}}\int_{|y|\ge\frac{x_{2n+2}-x_{2n+1}-2}{2}}e^{-\frac{|y_{2n+1}-y|^2}{4t}}\,dy.$$
Since $x_{n+1}-x_n\to+\infty$ as $n\to+\infty$, one gets that $\liminf_{x\to+\infty}v(t,x)=\alpha$ and $\limsup_{x\to+\infty}v(t,x)=\beta$ for each $t>0$ (and also for $t=0$ since $v(0,\cdot)=u_0$). Finally, one concludes that
$$\alpha_{\min}=\alpha\ \hbox{ and }\ \alpha_{\max}=\beta.$$
Theorem~\ref{mainTh} implies then that
\be\label{calphabeta*}
c_{\alpha}\le c_*\le c^*\le c_{\beta}.
\ee

\subsubsection*{First estimates of $u(t,ct+x)$ when $c\le c_{\alpha}$ or $c\ge c_{\beta}$}

Let us come back to the solution $u$ of (\ref{cauchy1d}) with initial condition~$u_0$ given by~(\ref{defu0}) and~(\ref{defu0bis}). The function~$u_0$ satisfies $\alpha\le u_0\le 1$ in $\R$, whence
$$\alpha\le u(t,\cdot)\le 1\ \hbox{ in }\R$$
for all $t\ge 0$ (since $f(\alpha)=f(1)=0$). Since $\liminf_{x\to+\infty}v(t,x)=\alpha$, $\limsup_{x\to+\infty}v(t,x)=\beta$ and
$$\liminf_{x\to+\infty}u(t,x)=\liminf_{x\to+\infty}v(t,x),\quad\limsup_{x\to+\infty}u(t,x)=\limsup_{x\to+\infty}v(t,x)$$
from Lemma~\ref{L2}, one concludes that
\be\label{infsup}
\inf_{\R}u(t,\cdot)=\liminf_{x\to+\infty}u(t,x)=\alpha<\beta=\limsup_{x\to+\infty}u(t,x)
\ee
for all $t\ge 0$. Furthermore, for each $\eta\in(0,1-\theta)$, there is a real number~$\xi$ such that $\underline{\varphi}_{\alpha-\eta,\eta}(x+\xi)\le u_0(x)$ for all $x\in\R$, under the notations~(\ref{underc}). Since $\underline{f}_{\eta}\le f$, the maximum principle implies that
$$\underline{\varphi}_{\alpha-\eta,\eta}(x-\underline{c}_{\alpha-\eta,\eta}t+\xi)\le u(t,x)\ \hbox{ for all }t\ge 0\hbox{ and }x\in\R,$$
whence $\liminf_{x\to-\infty}u(t,x)\ge\underline{\varphi}_{\alpha-\eta,\eta}(-\infty)=1-\eta$ for all $t\ge 0$. Since $\eta>0$ is arbitrarily small and $u\le 1$, one gets that
$$\forall\ t\ge 0,\quad\sup_{\R}u(t,\cdot)=\lim_{x\to-\infty}u(t,x)=1.$$\par
According to the notations~(\ref{cauchy1}) and~(\ref{cauchy2}), there holds
$$\alpha\le\underline{u}_{\alpha,0}(0,\cdot)\le u_0\le\overline{u}_{\beta}(0,\cdot)\le 1\ \hbox{ in }\R.$$
As a consequence,
\be\label{compar}
\alpha\le\underline{u}_{\alpha,0}(t,x)\le u(t,x)\le\overline{u}_{\beta}(t,x)\le 1\ \hbox{ for all }(t,x)\in[0,+\infty)\times\R.
\ee
Since
\be\label{x0}\baa{l}
\displaystyle{\mathop{\lim}_{t\to+\infty}}\Big(\displaystyle{\mathop{\sup}_{x\in\R}}\left|\underline{u}_{\alpha,0}(t,x)-\varphi_{\alpha}(x-c_{\alpha}t+\underline{x}_{\alpha,0})\right|\Big)\vspace{5pt}\\
\qquad\qquad\qquad\qquad\qquad\qquad=\displaystyle{\mathop{\lim}_{t\to+\infty}}\Big(\displaystyle{\mathop{\sup}_{x\in\R}}\left|\overline{u}_{\beta}(t,x)-\varphi_{\beta}(x-c_{\beta}t+\overline{x}_{\beta})\right|\Big)=0\eaa
\ee
and since all functions $\varphi_{\gamma}$ are continuous and decreasing in $\R$, it follows that
\be\label{estimates}\left\{\baa{rl}
\forall\ c<c_{\alpha},\ \forall\ A\in\R, & u(t,ct+\cdot)\displaystyle{\mathop{\longrightarrow}_{t\to+\infty}}1\ \hbox{ uniformly in }(-\infty,A],\vspace{5pt}\\
\forall\ A\in\R, & 
\alpha<\varphi_{\alpha}(A+\underline{x}_{\alpha,0})\le\displaystyle{\mathop{\liminf}_{t\to+\infty}}\Big(\displaystyle{\mathop{\inf}_{x\in(-\infty,A]}}u(t,c_{\alpha}t+x)\Big)\ (\le 1),\vspace{5pt}\\
\forall\ A\in\R, & (\alpha\le)\ \displaystyle{\mathop{\limsup}_{t\to+\infty}}\Big( \displaystyle{\mathop{\sup}_{x\in[A,+\infty)}}u(t,c_{\beta}t+x)\Big)\le\varphi_{\beta}(A+\overline{x}_{\beta})<1,\vspace{5pt}\\
\forall\ c>c_{\beta},\ \forall\ A\in\R, & (\alpha\le)\ \displaystyle{\mathop{\limsup}_{t\to+\infty}}\Big( \displaystyle{\mathop{\sup}_{x\in[A,+\infty)}}u(t,ct+x)\Big)\le\beta.\eaa\right.
\ee
Notice also that, for all $x\in\R$,
\be\label{alphax}
\alpha<\alpha_x:=\displaystyle{\mathop{\liminf}_{t\to+\infty}}\ u(t,c_{\alpha}t+x)\le\displaystyle{\mathop{\limsup}_{t\to+\infty}}\ u(t,c_{\alpha}t+x)\le1
\ee
and
\be\label{betax}
\alpha\le\displaystyle{\mathop{\liminf}_{t\to+\infty}}\ u(t,c_{\beta}t+x)\le\displaystyle{\mathop{\limsup}_{t\to+\infty}}\ u(t,c_{\beta}t+x)=:\beta_x<1.
\ee

\subsubsection*{Definition of the functions $t\mapsto\xi(t)$ and $x\mapsto\tau(x)$}

The function $u_0$ is Lipschitz-continuous, piecewise $C^1$, and the value
$$\xi_0=\frac{2(1-\theta)}{1-\alpha}\in(0,2)$$
is the unique real number such that $u_0(\xi_0)=\theta$. Furthermore, $u'_0(\xi_0)=-(1-\alpha)/2<0$. Remember also that, for each $t>0$, the function $u(t,\cdot)$ is continuous and $u(t,-\infty)=1$, $\limsup_{x\to+\infty}u(t,x)\le\beta<\theta$. Since the number of intersection points of the function~$u(t,\cdot)$ with the constant~$\theta$ (which is a solution of the same parabolic equation as $u$) is nonincreasing in time, one concludes that, for each $t\ge 0$, there is a unique $\xi(t)\in\R$ such that
$$u(t,\xi(t))=\theta,\quad u(t,\cdot)>\theta\ \hbox{ in }(-\infty,\xi(t)),\quad u(t,\cdot)<\theta\ \hbox{ in }(\xi(t),+\infty),$$
and $u_x(t,\xi(t))<0$ (with these notations, there holds $\xi(0)=\xi_0$). It follows from the implicit function theorem that~$\xi$ is a~$C^1$ function of~$t$. Lastly, from~(\ref{compar}) and~(\ref{x0}), there exists a constant $M\ge 0$ such that
\be\label{xit}
\forall\ t\ge 0,\quad c_{\alpha}t-M\le\xi(t)\le c_{\beta}t+M.
\ee\par
Therefore, for each $x\ge\xi_0$, the real number
$$\tau(x)=\min\ \{t\ge 0,\ \xi(t)=x\}$$
is well-defined. Notice that $\tau(x)>0$ for all $x>\xi_0$. For all $x\ge\xi_0$, there holds $\xi(\tau(x))=x$ and $u(\tau(x),x)=\theta$. Furthermore, for all $t\in[0,\tau(x)]$, one has $\xi(t)\le x$. As a consequence,
\be\label{theta}
\forall\ x\ge\xi_0,\ \forall\ t\in[0,\tau(x)],\quad u(t,\cdot)\le\theta\ \hbox{ in }[x,+\infty).
\ee\par
For any $\xi_0\le x_1<x_2$, there holds $u(\tau(x_2),x_1)>\theta$ since $x_1<x_2=\xi(\tau(x_2))$. But $u(0,x_1)\le\theta$ since $x_1\ge\xi_0$. Consequently, $\tau(x_1)<\tau(x_2)$. Thus, the function $\tau:[\xi_0,+\infty)\to[0,+\infty)$ is increasing.\par
Lastly, notice from (\ref{xit}) (applied at $t=\tau(x)$) implies that
\be\label{taux}
\forall\ x\ge\xi_0,\quad\frac{x-M}{c_{\beta}}\le\tau(x)\le\frac{x+M}{c_{\alpha}}.
\ee
In particular, $\lim_{n\to+\infty}\tau(x_n)=+\infty$, since $\lim_{n\to+\infty}x_n=+\infty$.

\subsubsection*{The key-lemma}

The key-point in the proof of Theorem~\ref{th1d} is the following lemma, the proof of which is postponed in the next subsection:

\begin{lem}\label{speeds}
Set $z_n=\sqrt{x_nx_{n+1}}$ for each $n\in\N$. For each $0<\epsilon<1$, there exists $\eta_0=\eta_0(\epsilon)>0$ such that the following holds: for all $\eta\in(0,\eta_0)$, there is $N=N(\epsilon,\eta)\in\N$ such that
$$\forall\ n\ge N,\quad\left\{\baa{l}
\left|\displaystyle{\frac{x_{2n}}{\tau(x_{2n})}}-c_{\beta}\right|+\left|\displaystyle{\frac{z_{2n}}{\tau(z_{2n})}}-c_{\alpha}\right|+\left|\displaystyle{\frac{x_{2n+1}}{\tau(x_{2n+1})}}-c_{\alpha}\right|+\left|\displaystyle{\frac{z_{2n+1}}{\tau(z_{2n+1})}}-c_{\beta}\right|\le\epsilon,\vspace{5pt}\\
u(\tau(x_{2n}),\cdot)\ge 1-\eta\ \hbox{ in }(-\infty,x_{2n}-\epsilon x_{2n}],\vspace{5pt}\\
u(\tau(x_{2n}),\cdot)\le \alpha+\eta\ \hbox{ in }[x_{2n}+\epsilon x_{2n},x_{2n+1}-\epsilon x_{2n+1}],\vspace{5pt}\\
u(\tau(z_{2n}),\cdot)\ge 1-\eta\ \hbox{ in }(-\infty,z_{2n}-\epsilon z_{2n}],\vspace{5pt}\\
u(\tau(z_{2n}),\cdot)\le \alpha+\eta\ \hbox{ in }[z_{2n}+\epsilon z_{2n},x_{2n+1}-\epsilon x_{2n+1}],\vspace{5pt}\\
u(\tau(x_{2n+1}),\cdot)\ge 1-\eta\ \hbox{ in }(-\infty,x_{2n+1}-\epsilon x_{2n+1}],\vspace{5pt}\\
|u(\tau(x_{2n+1}),\cdot)-\beta|\le\eta\ \hbox{ in }[x_{2n+1}+\epsilon x_{2n+1},x_{2n+2}-\epsilon x_{2n+2}],\vspace{5pt}\\
u(\tau(z_{2n+1}),\cdot)\ge 1-\eta\ \hbox{ in }(-\infty,z_{2n+1}-\epsilon z_{2n+1}],\vspace{5pt}\\
|u(\tau(z_{2n+1}),\cdot)-\beta|\le\eta\ \hbox{ in }[z_{2n+1}+\epsilon z_{2n+1},x_{2n+2}-\epsilon x_{2n+2}].\eaa\right.$$
\end{lem}

\subsubsection*{End of the proof of Theorem~\ref{th1d}}

First, let $c$ be any given speed such that $c<c_{\beta}$, let $x$ be any given real number and let us prove that $\limsup_{t\to+\infty}u(ct,t+x)=1$. Let $0<\epsilon<1$ be such that
$$c<(1-\epsilon)\times(c_{\beta}-\epsilon).$$
Let $\eta_0=\eta_0(\epsilon)$ be given by Lemma~\ref{speeds}. Pick any $\eta\in(0,\eta_0)$ and let $N=N(\epsilon,\eta)\in\N$ be given by Lemma~\ref{speeds}. Since $\tau(x_{2n})\to+\infty$ as $n\to+\infty$, there is $N_1\ge N$ such that
$$\forall\ n\ge N_1,\quad c+\displaystyle{\frac{x}{\tau(x_{2n})}}\le(1-\epsilon)\times(c_{\beta}-\epsilon).$$
For any $n\ge N_1$, it follows then from Lemma~\ref{speeds} that
$$c+\frac{x}{\tau(x_{2n})}\le(1-\epsilon)\times(c_{\beta}-\epsilon)\le(1-\epsilon)\times\frac{x_{2n}}{\tau(x_{2n})},$$
whence $c\ \!\tau(x_{2n})+x\le x_{2n}-\epsilon x_{2n}$. Thus,
$$\forall\ n\ge N_1,\quad u(\tau(x_{2n}),c\ \!\tau(x_{2n})+x)\ge1-\eta$$
from Lemma~\ref{speeds}. Since $\eta$ is arbitrary in $(0,\eta_0)$ and since $u\le 1$, one concludes that
\be\label{limsup1}
\forall\ c<c_{\beta},\ \forall\ x\in\R,\quad\limsup_{t\to+\infty}u(t,ct+x)=1.
\ee\par
Let now $c$ be any given speed such that $c>c_{\alpha}$, let $x$ be any given real number and let us prove that $\liminf_{t\to+\infty}u(ct,t+x)=\alpha$ and $\limsup_{t\to+\infty}u(t,ct+x)\ge\beta$. Let $0<\epsilon<1$ be such that
$$(1+\epsilon)\times(c_{\alpha}+\epsilon)<c,$$
let $\eta_0=\eta_0(\epsilon)$ be given by Lemma~\ref{speeds}, pick any $\eta\in(0,\eta_0)$ and let $N=N(\epsilon,\eta)\in\N$ be given by Lemma~\ref{speeds}. Since $\tau(z_{2n})\to+\infty$ and $\tau(x_{2n+1})\to+\infty$ as $n\to+\infty$, there is $N_1\ge N$ such that
$$(1+\epsilon)\times(c_{\alpha}+\epsilon)\le\min\left(c+\displaystyle{\frac{x}{\tau(z_{2n})}},c+\displaystyle{\frac{x}{\tau(x_{2n+1})}}\right).$$
Lemma~\ref{speeds} also implies that, for any $n\ge N_1$,
\be\label{z2n1}\left\{\baa{l}
(1+\epsilon)\times\displaystyle{\frac{z_{2n}}{\tau(z_{2n})}}\le(1+\epsilon)\times(c_{\alpha}+\epsilon)\le c+\displaystyle{\frac{x}{\tau(z_{2n})}},\vspace{5pt}\\
(1+\epsilon)\times\displaystyle{\frac{x_{2n+1}}{\tau(x_{2n+1})}}\le(1+\epsilon)\times(c_{\alpha}+\epsilon)\le c+\displaystyle{\frac{x}{\tau(x_{2n+1})}}.\eaa\right.
\ee
On the other hand,
$$\frac{x_{2n+1}}{\tau(z_{2n})}=\frac{z_{2n}}{\tau(z_{2n})}\times\frac{z_{2n}}{x_{2n}}=\frac{z_{2n}}{\tau(z_{2n})}\times\sqrt{\frac{x_{2n+1}}{x_{2n}}}\to+\infty\ \hbox{ as }n\to+\infty$$
from (\ref{xnn1}) and (\ref{taux}). Moreover,
$$\frac{x_{2n+2}}{\tau(x_{2n+1})}=\frac{x_{2n+1}}{\tau(x_{2n+1})}\times\frac{x_{2n+2}}{x_{2n+1}}\to+\infty\ \hbox{ as }n\to+\infty.$$
In particular, there exists $N_2\ge N_1$ such that
\be\label{z2n2}\forall\ n\ge N_2,\quad\left\{\baa{l}
c+\displaystyle{\frac{x}{\tau(z_{2n})}}\le(1-\epsilon)\times\displaystyle{\frac{x_{2n+1}}{\tau(z_{2n})}},\vspace{5pt}\\
c+\displaystyle{\frac{x}{\tau(x_{2n+1})}}\le(1-\epsilon)\times\displaystyle{\frac{x_{2n+2}}{\tau(x_{2n+1})}}.\eaa\right.
\ee
Eventually, it follows from (\ref{z2n1}) and (\ref{z2n2}) that
$$\forall\ n\ge N_2,\quad\left\{\baa{l}
z_{2n}+\epsilon z_{2n}\le c\ \!\tau(z_{2n})+x\le x_{2n+1}-\epsilon x_{2n+1},\vspace{5pt}\\
x_{2n+1}+\epsilon x_{2n+1}\le c\ \!\tau(x_{2n+1})+x\le x_{2n+2}-\epsilon x_{2n+2},\eaa\right.$$
whence
$$\forall\ n\ge N_2,\quad\left\{\baa{l}
u(\tau(z_{2n}),c\ \!\tau(z_{2n})+x)\le\alpha+\eta,\vspace{5pt}\\
|u(\tau(x_{2n+1}),c\ \!\tau(x_{2n+1})+x)-\beta|\le\eta\eaa\right.$$
from Lemma~\ref{speeds}. Since $\eta$ is arbitrary in $(0,\eta_0)$ and since $u\ge\alpha$, one concludes that
\be\label{liminfalpha}
\forall\ c>c_{\alpha},\ \forall\ x\in\R,\quad\liminf_{t\to+\infty}u(t,ct+x)=\alpha<\beta
\le\limsup_{t\to+\infty}u(t,ct+x).
\ee\par
Since the function $u$ is continuous, properties (\ref{limsup1}) and (\ref{liminfalpha}) yield:
\be\label{adherence}
\forall\ c\in(c_{\alpha},c_{\beta}),\ \forall\ x\in\R,\quad\left\{\lim_{t_k\to+\infty}u(t_k,ct_k+x)\right\}=[\alpha,1].
\ee
Notice that, from (\ref{calphabeta*}) and the general definitions of $c_*$ and $c^*$ given in the introduction, formula~(\ref{adherence}) implies in particular that
$$c_*=c_{\alpha}<c_{\beta}=c^*.$$
The second and third assertions in (\ref{estimates}) then yield (\ref{c**}).\par
Furthermore, property~(\ref{limsup1}) also implies that, for all $x\in\R$,
$$\left\{\lim_{t_k\to+\infty}u(t_k,c_{\alpha}t_k+x)\right\}=[\alpha_x,1],$$
where $\alpha_x=\liminf_{t\to+\infty}u(t,c_{\alpha}t+x)\in(\alpha,1]$ (see~(\ref{alphax})). Similarly, property (\ref{liminfalpha}) implies that, for all $x\in\R$, the real number $\beta_x\in[\alpha,1)$ given by (\ref{betax}), namely $\beta_x=\limsup_{t\to+\infty}u(t,c_{\beta}t+x)$, is such that $\beta_x\in[\beta,1)$ and
$$\left\{\lim_{t_k\to+\infty}u(t_k,c_{\beta}t_k+x)\right\}=[\alpha,\beta_x].$$\par
Lastly, for any speed $c>c_{\beta}$ and for any real number~$x$, it follows from the last assertion in (\ref{estimates}) and from (\ref{liminfalpha}) that
$$\forall\ c>c_{\beta},\ \forall\ x\in\R,\quad
\left\{\lim_{t_k\to+\infty}u(t_k,ct_k+x)\right\}=[\alpha,\beta].$$
Furthermore, (\ref{infsup}) and (\ref{estimates}) imply that
$$\forall\ c>c_{\beta},\ \forall\ A\in\R,\quad\displaystyle{\mathop{\lim}_{t\to+\infty}}\Big(\displaystyle{\mathop{\sup}_{x\in[A,+\infty)}}u(t,ct+x)\Big)=\beta.$$
That completes the proof of Theorem~\ref{th1d}.\hfill$\Box$

\begin{rem}\label{remspeeds}{\rm It follows from (\ref{xit}), (\ref{taux}) and Lemma~\ref{speeds} that
$$c_{\alpha}=\liminf_{t\to+\infty}\frac{\xi(t)}{t}<\limsup_{t\to+\infty}\frac{\xi(t)}{t}=c_{\beta}$$
and
$$\frac{1}{c_{\beta}}=\liminf_{x\to+\infty}\frac{\tau(x)}{x}<\limsup_{x\to+\infty}\frac{\tau(x)}{x}=\frac{1}{c_{\alpha}}.$$
In particular, there is no speed $c$ such that the function $t\mapsto\xi(t)-ct$ is bounded and there are no $\gamma<\theta$ and $x_0\in\R$ such that $u(t,\xi(t)+\cdot)$ converges as $t\to+\infty$ to a front~$\varphi_{\gamma}(\cdot+x_0)$. These properties are very different from the usual results of the literature, which are concerned with initial conditions $u_0$ converging to a constant as $x\to+\infty$.}
\end{rem}

%%%%%%%%%%%%%%%%%%%%%%%%%%%%%%%%%%%%%%%%

\subsection{Proof of Lemma~\ref{speeds}}

\subsubsection*{Choices of $\eta_0=\eta_0(\epsilon)$ and parameters depending on $\eta\in(0,\eta_0)$}

Let $0<\epsilon<1$ be given. Let $\rho>0$ be chosen so that
\be\label{defrho}
0<\rho<\frac{1}{2}\ \hbox{ and }\ c_{\beta}-\frac{\epsilon}{4}<\left(c_{\beta}^{-1}+\rho\ \!\epsilon\ \!c_{\alpha}^{-1}\right)^{-1}.
\ee
From (\ref{monospeeds}) and (\ref{cgammalim}), there exists $\eta_0=\eta_0(\epsilon)\in(0,\min(1-\theta,\theta-\beta))$ such that, for all $\eta\in(0,\eta_0)$,
\be\label{calphabeta}
\forall\ \eta\in(0,\eta_0),\ \left\{\baa{l}
0<c_{\alpha}\le\overline{c}_{\alpha+\eta,\eta}<c_{\alpha}+\displaystyle{\frac{\epsilon}{4}},\vspace{5pt}\\
c_{\beta}-\displaystyle{\frac{\epsilon}{4}}<\underline{c}_{\beta-\eta,\eta/2}\le c_{\beta}\le\overline{c}_{\beta+\eta,\eta}< c_{\beta}+\displaystyle{\frac{\epsilon}{4}},\vspace{5pt}\\
1-\displaystyle{\frac{c_{\alpha}}{\overline{c}_{\alpha+\eta,\eta}}}\le
\displaystyle{\frac{\epsilon}{8}},\vspace{5pt}\\
1-\displaystyle{\frac{\underline{c}_{\beta-\eta,\eta/2}}{\overline{c}_{\beta+\eta,\eta}}}\le\displaystyle{\frac{\rho\ \!\epsilon}{4}}\le\displaystyle{\frac{\epsilon}{8}},\vspace{5pt}\\
c_{\beta}-\displaystyle{\frac{\epsilon}{4}}<\left(\underline{c}_{\beta-\eta,\eta/2}^{-1}+\rho\ \!\epsilon\ \!\underline{c}_{\alpha-\eta,3\eta/4}^{-1}\right)^{-1}.\eaa\right.
\ee\par
In the sequel, let $\eta$ be any given real number in the interval $(0,\eta_0)$. Remember that the pairs $(\underline{c}_{\alpha-\eta,\eta/4},\underline{\varphi}_{\alpha-\eta,\eta/4})$, $(\underline{c}_{\alpha-\eta,3\eta/4},\underline{\varphi}_{\alpha-\eta,3\eta/4})$ and $(\underline{c}_{\beta-\eta,\eta/2},\underline{\varphi}_{\beta-\eta,\eta/2})$ solve~(\ref{underc}) with nonlinearities $\underline{f}_{\eta/4}$, $\underline{f}_{3\eta/4}$ and $\underline{f}_{\eta/2}$ respectively, and limit values
$$\left\{\baa{l}
\underline{\varphi}_{\alpha-\eta,\eta/4}(-\infty)=1-\displaystyle{\frac{\eta}{4}}>\alpha-\eta=\underline{\varphi}_{\alpha-\eta,\eta/4}(+\infty),\vspace{5pt}\\
\underline{\varphi}_{\alpha-\eta,3\eta/4}(-\infty)=1-\displaystyle{\frac{3\eta}{4}}>\alpha-\eta=\underline{\varphi}_{\alpha-\eta,3\eta/4}(+\infty),\vspace{5pt}\\
\underline{\varphi}_{\beta-\eta,\eta/2}(-\infty)=1-\displaystyle{\frac{\eta}{2}}>\beta-\eta=\underline{\varphi}_{\beta-\eta,\eta/2}(+\infty)\eaa\right.$$
and that the pairs 
$(\overline{c}_{\alpha+\eta,\eta},\overline{\varphi}_{\alpha+\eta,\eta})$ and $(\overline{c}_{\beta+\eta,\eta},\overline{\varphi}_{\beta+\eta,\eta})$ solve~(\ref{overc}) with nonlinearity $\overline{f}_{\eta}$ and limit values
$$\left\{\baa{l}
\overline{\varphi}_{\alpha+\eta,\eta}(-\infty)=1+\eta>\alpha+\eta=\overline{\varphi}_{\alpha+\eta,\eta}(+\infty),\vspace{5pt}\\
\overline{\varphi}_{\beta+\eta,\eta}(-\infty)=1+\eta>\beta+\eta=\overline{\varphi}_{\beta+\eta,\eta}(+\infty).\eaa\right.$$
There exists a real number $A=A(\eta)>0$, which is fixed in the sequel, such that
\be\label{choiceA}\left\{\baa{rcll}
\overline{\varphi}_{\alpha+\eta,\eta} & \ge &  1 & \hbox{ in }(-\infty,-A],\vspace{5pt}\\
\varphi_{\alpha} & \ge & 1-\displaystyle{\frac{\eta}{8}}\ \ge\ \theta+\displaystyle{\frac{\eta}{8}} & \hbox{ in }(-\infty,-A],\vspace{5pt}\\
\overline{\varphi}_{\beta+\eta,\eta} & \ge & 1 & \hbox{ in }(-\infty,-A],\vspace{5pt}\\
\underline{\varphi}_{\alpha-\eta,\eta/4} & \le & \alpha & \hbox{ in }[A,+\infty),\vspace{5pt}\\
\underline{\varphi}_{\alpha-\eta,\eta/4} & \ge & 1-\displaystyle{\frac{\eta}{2}} & \hbox{ in }(-\infty,-A],\vspace{5pt}\\
\underline{\varphi}_{\beta-\eta,\eta/2} & \le & \beta-\displaystyle{\frac{\eta}{2}} & \hbox{ in }[A,+\infty),\vspace{5pt}\\
\underline{\varphi}_{\beta-\eta,\eta/2} & \ge & 1-\displaystyle{\frac{3\eta}{4}} & \hbox{ in }(-\infty,-A],\vspace{5pt}\\
\underline{\varphi}_{\alpha-\eta,3\eta/4} & \le & \alpha & \hbox{ in }[A,+\infty),\vspace{5pt}\\
\underline{\varphi}_{\alpha-\eta,3\eta/4} & \ge & 1-\eta & \hbox{ in }(-\infty,-A].\eaa\right.
\ee
Because of (\ref{conv1}), there exists also a time $T=T(\eta)\ge 0$ such that
\be\label{choiceT}
\forall\ t\ge T,\ \forall\ x\in\R,\quad|\underline{u}_{\alpha,\eta}(t,x)-\varphi_{\alpha}(x-c_{\alpha}t+\underline{x}_{\alpha,\eta})|\le\frac{\eta}{8}.
\ee

\subsubsection*{Comparisons with solutions of heat equations}

Let $v$ be the solution of the heat equation~(\ref{heatv}) with initial condition~$u_0$. We know that $\alpha\le u,v\le 1$ in $[0,+\infty)\times\R$. Furthermore, since $f\ge 0$, one gets that
$$\forall\ (t,x)\in[0,+\infty)\times\R,\quad 0\le\alpha\le v(t,x)\le u(t,x)\le 1.$$
On the other hand, for any given $x\in(\xi_0,+\infty)$, there holds $u_t(t,y)=u_{yy}(t,y)$ for all $(t,y)\in(0,\tau(x)]\times[x,+\infty)$ and $u(t,x)\le\theta$ for all $t\in[0,\tau(x)]$ from~(\ref{theta}). The maximum principle implies that
$$\forall\ (t,y)\in[0,\tau(x)]\times[x,+\infty),\quad u(t,y)\le v(t,y)+w(t,y),$$
where $w$ solves the heat equation $w_t=w_{yy}$ in $(0,+\infty)\times(x,+\infty)$, with $w(0,y)=0$ in~$(x,+\infty)$ and $w(t,x)=\theta$ for all $t>0$. The function $w$ is explicitely given by
$$\forall\ (t,y)\in(0,+\infty)\times[x,+\infty),\quad w(t,y)=\frac{2\theta}{\sqrt{\pi}}\int_{\frac{y-x}{2\sqrt{t}}}^{+\infty}e^{-z^2}dz.$$
Finally,
$$\forall\ x>\xi_0,\ \forall\ (t,y)\in(0,\tau(x)]\times[x,+\infty),\quad u(t,y)\le v(t,y)+\frac{2\theta}{\sqrt{\pi}}\int_{\frac{y-x}{2\sqrt{t}}}^{+\infty}e^{-z^2}dz.$$
Let now $B=B(\eta)>0$ be given so that
\be\label{choiceB}
\frac{1}{\sqrt{\pi}}\int_B^{+\infty}e^{-z^2}dz\le\frac{\eta}{4},
\ee
and $\xi_1=\xi_1(\epsilon,\eta)>\xi_0\ (>0)$ be such that
\be\label{choicexi1}
\forall\ x\ge\xi_1,\quad\min\left(\frac{x^{3/4}-1}{2\sqrt{\tau(x+x^{3/4})}},\frac{x-2}{2\sqrt{\tau(x)}}\right)\ge B.
\ee
The choice of $\xi_1$ is possible because of (\ref{taux}). In particular, there holds $x^{3/4}/(2\sqrt{\tau(x)})\ge B$ for all $x\ge\xi_1$, since $\tau$ is increasing. Thus,
\be\label{xi1}
\forall\ x\ge\xi_1,\ \forall\ t\in[0,\tau(x)],\quad u(t,\cdot)\le v(t,\cdot)+\frac{2\ \!\theta\ \!\eta}{4}\le v(t,\cdot)+\frac{\eta}{2}\ \hbox{ in }[x+x^{3/4},+\infty).
\ee
Notice indeed that the above inequality is immediate at time $t=0$.\par
Furthermore, for all $x\ge\max(\xi_1,2)$ and $(t,y)\in(0,\tau(x)]\times[x,+\infty)$, there holds that
$$\baa{rcl}
v(t,y)\!-\!\beta & \!\!=\!\! & \displaystyle{\frac{1}{\sqrt{4\pi t}}}\displaystyle{\int_{-\infty}^{+\infty}}e^{-\frac{|y-z|^2}{4t}}(u_0(z)-\beta)dz\vspace{5pt}\\
& \!\!\le\!\! & \displaystyle{\frac{1}{\sqrt{4\pi t}}}\displaystyle{\int_{-\infty}^2}e^{-\frac{|y-z|^2}{4t}}dz=\displaystyle{\frac{1}{\sqrt{\pi}}}\displaystyle{\int_{\frac{y-2}{2\sqrt{t}}}^{+\infty}}e^{-z^2}dz\le\displaystyle{\frac{1}{\sqrt{\pi}}}\displaystyle{\int_{\frac{x-2}{2\sqrt{\tau(x)}}}^{+\infty}}e^{-z^2}dz\le\displaystyle{\frac{\eta}{4}}\eaa$$
from (\ref{choiceB}) and (\ref{choicexi1}). Hence, it follows from (\ref{xi1}) that
\be\label{betaeta}
\forall\ x\ge\max(\xi_1,2),\ \forall\ t\in[0,\tau(x)],\ \  u(t,\cdot)\le\beta+\frac{\eta}{4}+\frac{\eta}{2}\le\beta+\eta\ \hbox{ in }[x+x^{3/4},+\infty),
\ee
where the above inequality also holds immediately at time $t=0$.

\subsubsection*{Choice of a first iteration point $x_{2N_0}$}

Remember that $z_m=\sqrt{x_mx_{m+1}}$ for each $m\in\N$, and that
$$\lim_{m\to+\infty}x_m=\lim_{m\to+\infty}\frac{x_{m+1}}{x_m}=\lim_{m\to+\infty}\frac{z_m}{x_m}=\lim_{m\to+\infty}\frac{x_{m+1}}{z_m}=+\infty.$$
Let $N_0=N_0(\epsilon,\eta)\in\N$ be such that
\be\label{choiceN0}\forall\ \!m\!\ge\! 2\ \!N_0,\!\left\{\baa{l}
\!\xi_1\!+\!2\!+\!A\!+\!M\!+\!\overline{c}_{\alpha+\eta,\eta}T\!+\!|\underline{x}_{\alpha,\eta}|\!+\!\displaystyle{\frac{8\ \!\!A}{\epsilon}}\!+\!A^{4/3}\!+\!\displaystyle{\frac{1}{(\rho\epsilon)^4}}\!+\!\displaystyle{\frac{4\ \!A}{(1\!-\!2\rho)\epsilon}}\le x_m,\vspace{5pt}\\
\!x_m<x_m+x_m^{3/4}\le x_m+\epsilon x_m\le 4\ \!x_m\le z_m\vspace{5pt}\\
\!z_m+z_m^{3/4}\le z_m+\epsilon z_m\le x_{m+1}-\epsilon x_{m+1}\le x_{m+1}-2\ \!x_{m+1}^{3/4}<x_{m+1},\vspace{5pt}\\
\!4\ \!x_m\le\left(6+\displaystyle{\frac{3\ \!c_{\beta}}{c_{\alpha}}}\right)x_m\le\displaystyle{\frac{\epsilon z_m}{4}}\le\displaystyle{\frac{\rho\epsilon x_{m+1}}{2}},\vspace{5pt}\\
\!4\ \!x_m+2\ \!x_{m+1}^{3/4}\le\displaystyle{\frac{\epsilon x_{m+1}}{4}},\vspace{5pt}\\
\!A\le\displaystyle{\frac{\rho\epsilon x_m}{4}}-2\ \!x_m^{3/4},\eaa\right.
\ee
where $\underline{x}_{\alpha,\eta}\in\R$, $M\ge 0$, $\rho\in(0,1/2)$, $A\ge 0$ and $T\ge 0$ are given in (\ref{conv1}), (\ref{xit}), (\ref{defrho}), (\ref{choiceA}) and (\ref{choiceT}).

\subsubsection*{Estimates of $v$ in intervals of the type $[x+x^{3/4},x_{m+1}-x_{m+1}^{3/4}]$}

Choose any integer $n$ such that $n\ge N_0$, any real number $x\in[x_{2n},x_{2n+1}-x_{2n+1}^{3/4}]$, any real number $t\in(0,\tau(x+x^{3/4})]$ and any real number $y\in[x+x^{3/4},x_{2n+1}-x_{2n+1}^{3/4}]$ (when this interval is not empty). Since $v(0,\cdot)=u_0=\alpha$ in the interval $[x_{2n}+1,x_{2n+1}-1]$, there holds that
$$\baa{rcl}
|v(t,y)-\alpha| & \!\le\! & \displaystyle{\frac{1}{\sqrt{4\pi t}}}\times\left(\displaystyle{\int_{-\infty}^{x_{2n}+1}}e^{-\frac{|y-z|^2}{4t}}|u_0(z)-\alpha|\,dz+\displaystyle{\int_{x_{2n+1}-1}^{+\infty}}e^{-\frac{|y-z|^2}{4t}}|u_0(z)-\alpha|\,dz\right)\vspace{5pt}\\
& \!\le\! & \displaystyle{\frac{1}{\sqrt{\pi}}}\times\left(\displaystyle{\int_{-\infty}^{\frac{x_{2n}+1-y}{2\sqrt{t}}}}e^{-z^2}dz+\displaystyle{\int_{\frac{x_{2n+1}-1-y}{2\sqrt{t}}}^{+\infty}}e^{-z^2}\,dz\right)\!,\eaa$$
while
$$\frac{x_{2n}\!+\!1\!-\!y}{2\sqrt{t}}\le\frac{-x^{3/4}\!+\!1}{2\sqrt{\tau(x+x^{3/4})}}\le-B\le B\le\frac{x^{3/4}\!-\!1}{2\sqrt{\tau(x\!+\!x^{3/4})}}\le\frac{x_{2n+1}^{3/4}\!-\!1}{2\sqrt{t}}\le\frac{x_{2n+1}\!-\!1\!-\!y}{2\sqrt{t}}$$
from (\ref{choicexi1}) and (\ref{choiceN0}). It follows then from (\ref{choiceB}) that
\be\label{vpair}\baa{l}
\forall\ n\ge N_0,\ \forall\ x\in[x_{2n},x_{2n+1}-x_{2n+1}^{3/4}],\ \forall\ t\in[0,\tau(x+x^{3/4})],\vspace{5pt}\\
\qquad\qquad\qquad\qquad\qquad\qquad\qquad|v(t,\cdot)-\alpha|\le\displaystyle{\frac{\eta}{2}}\ \hbox{ in }[x+x^{3/4},x_{2n+1}-x_{2n+1}^{3/4}],\eaa
\ee
provided that the space interval in not empty. Similarly, since $v(0,\cdot)=u_0=\beta$ in the interval $[x_{2n+1}+1,x_{2n+2}-1]$, one gets that
\be\label{vimpair}\baa{l}
\forall\ n\ge N_0,\ \forall\ x\in[x_{2n+1},x_{2n+2}-x_{2n+2}^{3/4}],\ \forall\ t\in[0,\tau(x+x^{3/4})],\vspace{5pt}\\
\qquad\qquad\qquad\qquad\qquad\qquad\qquad|v(t,\cdot)-\beta|\le\displaystyle{\frac{\eta}{2}}\ \hbox{ in }[x+x^{3/4},x_{2n+2}-x_{2n+2}^{3/4}],\eaa
\ee
provided that the space interval in not empty.

\subsubsection*{Refined estimates of $u(\tau(x),\cdot)$ in intervals of the type $[x+x^{3/4},x_{m+1}-x_{m+1}^{3/4}]$}

Let $n\ge N_0$ be given. Let us first show that $\alpha\le u(t,\cdot)\le\alpha+\eta$ in $[x_{2n}+x_{2n}^{3/4},x_{2n+1}-x_{2n+1}^{3/4}]$ for all $t\in(0,\tau(x_{2n})]$. This would then, in particular, yield the same inequality, at time $t=\tau(x_{2n})$, in the smaller interval $[x_{2n}+\epsilon x_{2n},x_{2n+1}-\epsilon x_{2n+1}]$, from the choice of $N_0$ in~(\ref{choiceN0}). Remember that the lower bound $u(t,x)\ge\alpha$ always holds. Furthermore, since $x_{2n}\ge\xi_1$ and $\tau(x_{2n})\le\tau(x_{2n}+x_{2n}^{3/4})$, properties~(\ref{xi1}) and (\ref{vpair}) --with $x=x_{2n}$-- imply that
$$\ \forall\ t\in[0,\tau(x_{2n})],\ \forall\ y\in[x_{2n}+x_{2n}^{3/4},x_{2n+1}-x_{2n+1}^{3/4}],\quad u(t,y)\le v(t,y)+\frac{\eta}{2}\le\alpha+\eta.$$
Eventually,
\be\label{alphaeta}
\forall\ n\ge N_0,\ \forall\ t\in[0,\tau(x_{2n})],\quad\alpha\le u(t,\cdot)\le\alpha+\eta\ \hbox{ in }[x_{2n}+x_{2n}^{3/4},x_{2n+1}-x_{2n+1}^{3/4}].
\ee\par
With the same arguments, the following estimates hold:
\be\label{right}\left\{\baa{ll}
\alpha\le u\le\alpha+\eta & \hbox{ in }[0,\tau(z_{2n})]\times[z_{2n}+z_{2n}^{3/4},x_{2n+1}-x_{2n+1}^{3/4}],\vspace{5pt}\\
\beta-\displaystyle{\frac{\eta}{2}}\le v\le u\le\beta+\eta & \hbox{ in }[0,\tau(x_{2n+1})]\times[x_{2n+1}+x_{2n+1}^{3/4},x_{2n+2}-x_{2n+2}^{3/4}],\vspace{5pt}\\
\beta-\displaystyle{\frac{\eta}{2}}\le v\le u\le\beta+\eta & \hbox{ in }[0,\tau(z_{2n+1})]\times[z_{2n+1}+z_{2n+1}^{3/4},x_{2n+2}-x_{2n+2}^{3/4}]\eaa\right.
\ee
for all $n\ge N_0$. The last two properties follow from (\ref{xi1}) and (\ref{vimpair}) applied with $x=x_{2n+1}$ and $x=z_{2n+1}$ respectively. Notice that these three properties then hold a fortiori in the smaller space intervals $[z_{2n}+\epsilon z_{2n},x_{2n+1}-\epsilon x_{2n+1}]$, $[x_{2n+1}+\epsilon x_{2n+1},x_{2n+2}-\epsilon x_{2n+2}]$ and $[z_{2n+1}+\epsilon z_{2n+1},x_{2n+2}-\epsilon x_{2n+2}]$ respectively. Actually, one gets more generally that
\be\label{upper}
\forall\ \!x\!\in\![x_{2n},x_{2n+1}\!-\!x_{2n+1}^{3/4}],\ \alpha\le u\le\alpha\!+\!\eta\hbox{ in }[0,\tau(x)]\!\times\![x\!+\!x^{3/4},x_{2n+1}\!-\!x_{2n+1}^{3/4}]
\ee
and
$$\forall\ \!x\!\in\![x_{2n+1},x_{2n+2}\!-\!x_{2n+2}^{3/4}],\ \beta\!-\!\displaystyle{\frac{\eta}{2}}\le u\le\beta\!+\!\eta\hbox{ in }[0,\tau(x)]\!\times\![x\!+\!x^{3/4},x_{2n+2}\!-\!x_{2n+2}^{3/4}]$$
for all $n\ge N_0$, provided that the space intervals are not empty.

\subsubsection*{From time $t=\tau(x_{2n})$ to time $t=\tau(x_{2n+1})$}

The heart of the proof of Lemma~\ref{speeds} consists in estimating from below $u(\tau(x),\cdot)$ on~$(-\infty,x-\epsilon x]$ and estimating some ratios $x/\tau(x)$, for $x=x_{2n}$,~$z_{2n}$,~$x_{2n+1}$ and~$z_{2n+1}$. We will do that by induction on~$n$ and step by step, from time $\tau(x_{2n})$ to time $\tau(x_{2n+1})$, and from time $\tau(x_{2n+1})$ to time~$\tau(x_{2n+2})$.\hfill\break

\noindent{\underline{Step 1: lower bound of $\tau(x)$ for $x\in[x_{2n},x_{2n+1}-2\ \!x_{2n+1}^{3/4}]$}}. Choose any integer $n$ such that $n\ge N_0$. There holds
$$u(\tau(x_{2n}),\cdot)\le\alpha+\eta\ \hbox{ in }[x_{2n}+x_{2n}^{3/4},x_{2n+1}-x_{2n+1}^{3/4}]$$
from~(\ref{alphaeta}). Moreover, $u(\tau(x_{2n}),\cdot)\le 1$ in $\R$. It follows then from the first assertion in~(\ref{choiceA}) and from the inequality $\overline{\varphi}_{\alpha+\eta,\eta}\ge\alpha+\eta$ in $\R$ that
$$u(\tau(x_{2n}),x)\le\overline{\varphi}_{\alpha+\eta,\eta}(x-(x_{2n}+x_{2n}^{3/4})-A)\ \hbox{ for all }x\in(-\infty,x_{2n+1}-x_{2n+1}^{3/4}].$$
Furthermore, since
$$x_{2n}\le x_{2n+1}-2\ \!x_{2n+1}^{3/4}\ \hbox{ and }\ (x_{2n+1}-2\ \!x_{2n+1}^{3/4})+(x_{2n+1}-2\ \!x_{2n+1}^{3/4})^{3/4}\le x_{2n+1}-x_{2n+1}^{3/4},$$
it follows from (\ref{upper}), applied at $x=x_{2n+1}-2\ \!x_{2n+1}^{3/4}$, that
$$\forall\ t\in[0,\tau(x_{2n+1}-2\ \!x_{2n+1}^{3/4})],\quad u(t,x_{2n+1}-x_{2n+1}^{3/4})\le\alpha+\eta.$$
Since $\overline{f}_{\eta}\ge f$, the function $\overline{\varphi}_{\alpha+\eta,\eta}(x-\overline{c}_{\alpha+\eta,\eta}t)$ is a supersolution of the equation satisfied by $u$. Since $\overline{\varphi}_{\alpha+\eta,\eta}\ge\alpha+\eta$ in $\R$, the maximum principle applied in the set where $(t,x)\in[\tau(x_{2n}),\tau(x_{2n+1}-2\ \!x_{2n+1}^{3/4})]\times(-\infty,x_{2n+1}-x_{2n+1}^{3/4}]$ then yields
$$u(t,x)\le\overline{\varphi}_{\alpha+\eta,\eta}\left(x-x_{2n}-x_{2n}^{3/4}-A-\overline{c}_{\alpha+\eta,\eta}(t-\tau(x_{2n}))\right).$$
for all $(t,x)\in[\tau(x_{2n}),\tau(x_{2n+1}-2\ \!x_{2n+1}^{3/4})]\times(-\infty,x_{2n+1}-x_{2n+1}^{3/4}]$. In particular, by choosing $t=\tau(x)$ and $x\in[x_{2n},x_{2n+1}-2\ \!x_{2n+1}^{3/4}]$, one has
$$\theta=u(\tau(x),x)\le\overline{\varphi}_{\alpha+\eta,\eta}\left(x-x_{2n}-x_{2n}^{3/4}-A-\overline{c}_{\alpha+\eta,\eta}(\tau(x)-\tau(x_{2n}))\right)\!.$$
But $\overline{\varphi}_{\alpha+\eta,\eta}$ is decreasing and equals $\theta$ at $0$. Hence,
\be\label{tauinf}\baa{rcl}
\forall\ x\in[x_{2n},x_{2n+1}-2\ \!x_{2n+1}^{3/4}],\quad\tau(x) & \ge & \displaystyle{\frac{x-x_{2n}-x_{2n}^{3/4}-A}{\overline{c}_{\alpha+\eta,\eta}}}+\tau(x_{2n})\vspace{5pt}\\
& \ge & \displaystyle{\frac{x-3\ \!x_{2n}}{\overline{c}_{\alpha+\eta,\eta}}}+\tau(x_{2n}),\eaa
\ee
since $x_{2n}\ge\max(1,A)$ from (\ref{choiceN0}).\hfill\break

\noindent{\underline{Step 2: upper bound of $\tau(x)$ for $x\ge 3\ \!x_{2n}+\overline{c}_{\alpha+\eta,\eta}T$}}. Let $n$ be any given integer such that $n\ge N_0$, and let $X\ge 0$ be such that
\be\label{hyprec}
u(\tau(x_{2n}),\cdot)\ge 1-\eta\ \hbox{ in }(-\infty,x_{2n}-\epsilon'x_{2n}-X],
\ee
where we set
$$\epsilon'=\frac{\epsilon}{2}.$$
Notice that such a $X\ge 0$ always exists since $u(\tau(x_{2n}),-\infty)=1$. Owing to the definition of $\underline{u}_{\alpha,\eta}$ in (\ref{cauchy1}), and since $u(\tau(x_{2n}),\cdot)\ge\alpha$ in $\R$, there holds then
$$\forall\ x\in\R,\quad\underline{u}_{\alpha,\eta}\left(0,x-(x_{2n}-\epsilon'x_{2n}-X)+2\right)\le u(\tau(x_{2n}),x).$$
But $\underline{u}_{\alpha,\eta}$ is a subsolution of the equation satisfied by $u$, since $\underline{f}_{\eta}\le f$. Thus,
$$\forall\ t\ge\tau(x_{2n}),\ \forall\ x\in\R,\ \ \underline{u}_{\alpha,\eta}\left(t-\tau(x_{2n}),x-x_{2n}+\epsilon'x_{2n}+X+2\right)\le u(t,x).$$
Hence, for all $(t,x)\in[T+\tau(x_{2n}),+\infty)\times\R$,
\be\label{varphiinf}
\varphi_{\alpha}\left(x-x_{2n}+\epsilon'x_{2n}+X+2+\underline{x}_{\alpha,\eta}-c_{\alpha}(t-\tau(x_{2n})\right)-\frac{\eta}{8}\le u(t,x)
\ee
from (\ref{choiceT}). Since $3\ \!x_{2n}+\overline{c}_{\alpha+\eta,\eta}T\le 4\ \!x_{2n}\le x_{2n+1}-2\ \!x_{2n+1}^{3/4}$ from (\ref{choiceN0}), and since $\tau$ is increasing, there holds
\be\label{compartau}
\forall\ x\ge 3\ \!x_{2n}+\overline{c}_{\alpha+\eta,\eta}T,\quad\tau(x)\ge\tau(3\ \!x_{2n}+\overline{c}_{\alpha+\eta,\eta}T)\ge T+\tau(x_{2n}),
\ee
where the last inequality follows from (\ref{tauinf}). In particular, by choosing $x\ge3\ \!x_{2n}+\overline{c}_{\alpha+\eta,\eta}T$ and $t=\tau(x)\ge T+\tau(x_{2n})$ in~(\ref{varphiinf}), one gets that
$$\varphi_{\alpha}\left(x-x_{2n}+\epsilon'x_{2n}+X+2+\underline{x}_{\alpha,\eta}-c_{\alpha}(\tau(x)-\tau(x_{2n})\right)-\frac{\eta}{8}\le u(\tau(x),x)=\theta,$$
whence
$$x-x_{2n}+\epsilon'x_{2n}+X+2+\underline{x}_{\alpha,\eta}-c_{\alpha}(\tau(x)-\tau(x_{2n})\ge-A$$
from the second assertion in (\ref{choiceA}) and since $\varphi_{\alpha}$ is decreasing. Thus,
\be\label{tausup}\baa{rcl}
\forall\ x\ge3\ x_{2n}+\overline{c}_{\alpha+\eta,\eta}T,\quad\tau(x) & \le & \displaystyle{\frac{x-x_{2n}+\epsilon'x_{2n}+X+A+2}{c_{\alpha}}}+\tau(x_{2n})\vspace{5pt}\\
& \le & \displaystyle{\frac{x+\epsilon'x_{2n}+X+A+2+M}{c_{\alpha}}}\vspace{5pt}\\
& \le & \displaystyle{\frac{x+2\ \!x_{2n}+X}{c_{\alpha}}}\eaa
\ee
from (\ref{taux}) and (\ref{choiceN0}), and since $\epsilon'=\epsilon/2<1/2<1$.\hfill\break

\noindent{\underline{Step 3: estimates of $\tau(z_{2n})$ and $\tau(x_{2n+1})$ and lower bound of $u$ on the left of $z_{2n}$ and $x_{2n+1}$}}. Notice that 
$$3\ \!x_{2n}+\overline{c}_{\alpha+\eta,\eta}T\le 4\ \!x_{2n}\le z_{2n}\le x_{2n+1}-2\ \!x_{2n+1}^{3/4}\le x_{2n+1}$$
because of (\ref{choiceN0}). As a consequence, it follows from (\ref{tauinf}), (\ref{tausup}) and the monotonicity of $\tau$, that
\be\label{tauznxn}\left\{\baa{l}
\displaystyle{\frac{z_{2n}-3\ \!x_{2n}}{\overline{c}_{\alpha+\eta,\eta}}}\le\displaystyle{\frac{z_{2n}-3\ \!x_{2n}}{\overline{c}_{\alpha+\eta,\eta}}}+\tau(x_{2n})\le \tau(z_{2n})\le\displaystyle{\frac{z_{2n}+2\ \!x_{2n}+X}{c_{\alpha}}},\vspace{5pt}\\
\displaystyle{\frac{x_{2n+1}-2\ \!x_{2n+1}^{3/4}-3\ \!x_{2n}}{\overline{c}_{\alpha+\eta,\eta}}}\le\displaystyle{\frac{x_{2n+1}-2\ \!x_{2n+1}^{3/4}-3\ \!x_{2n}}{\overline{c}_{\alpha+\eta,\eta}}}+\tau(x_{2n})\le\cdots\vspace{5pt}\\
\qquad\qquad\cdots\le\tau(x_{2n+1}-2\ \!x_{2n+1}^{3/4})\le
\tau(x_{2n+1})\le\displaystyle{\frac{x_{2n+1}+2\ \!x_{2n}+X}{c_{\alpha}}},\eaa\right.
\ee
provided that (\ref{hyprec}) holds. Since $x_{2n+1}\ge z_{2n}\ge 3\ \!x_{2n}+\overline{c}_{\alpha+\eta,\eta}T$, it follows from (\ref{varphiinf}) and~(\ref{compartau}) that
\be\label{uphi}\left\{\baa{l}
\!\!u(\tau(z_{2n}),x)\ge\varphi_{\alpha}\left(x\!-\!x_{2n}\!+\!\epsilon'x_{2n}\!+\!X\!+\!2\!+\!\underline{x}_{\alpha,\eta}\!-\!c_{\alpha}(\tau(z_{2n})\!-\!\tau(x_{2n}))\right)\!-\!\displaystyle{\frac{\eta}{8}},\vspace{5pt}\\
\!\!u(\tau(x_{2n+1}),x)\ge\varphi_{\alpha}\left(x\!-\!x_{2n}\!+\!\epsilon'x_{2n}\!+\!X\!+\!2\!+\!\underline{x}_{\alpha,\eta}\!-\!c_{\alpha}(\tau(x_{2n+1})\!-\!\tau(x_{2n}))\right)\!-\!\displaystyle{\frac{\eta}{8}}\eaa\right.
\ee
for all $x\in\R$. On the other hand, for all $x\le z_{2n}-\epsilon'z_{2n}-X=z_{2n}-\epsilon z_{2n}/2-X$, there holds
$$\baa{l}
x-x_{2n}+\epsilon'x_{2n}+X+2+\underline{x}_{\alpha,\eta}-c_{\alpha}(\tau(z_{2n})-\tau(x_{2n}))\vspace{5pt}\\
\qquad\qquad\qquad\qquad\qquad\qquad\qquad\le\left(1-\displaystyle{\frac{c_{\alpha}}{\overline{c}_{\alpha+\eta,\eta}}}-\displaystyle{\frac{\epsilon}{2}}\right)z_{2n}+3\ \!x_{2n}+2+\underline{x}_{\alpha,\eta}\vspace{5pt}\\
\qquad\qquad\qquad\qquad\qquad\qquad\qquad\le\left(\displaystyle{\frac{\epsilon}{8}}-\displaystyle{\frac{\epsilon}{2}}\right)z_{2n}+4\ \!x_{2n}\le-\displaystyle{\frac{\epsilon}{8}}\ \!z_{2n}\le-A
\eaa$$
from (\ref{calphabeta}), (\ref{choiceN0}), (\ref{tauznxn}) and since $-x_{2n}+\epsilon'x_{2n}\le 0$. Thus,
\be\label{tauz2n}
u(\tau(z_{2n}),\cdot)\ge\varphi_{\alpha}(-A)-\frac{\eta}{8}\ge1-\frac{\eta}{4}\ge 1-\eta\ \hbox{ in }(-\infty,z_{2n}-\epsilon'z_{2n}-X]
\ee
from (\ref{uphi}), from the second assertion in (\ref{choiceA}), and since $\varphi_{\alpha}$ is decreasing. Similarly, for all $x\le x_{2n+1}-\epsilon'x_{2n+1}-X$, there holds
$$\baa{l}
x-x_{2n}+\epsilon'x_{2n}+X+2+\underline{x}_{\alpha,\eta}-c_{\alpha}(\tau(x_{2n+1})-\tau(x_{2n}))\vspace{5pt}\\
\qquad\qquad\qquad\qquad\le\left(1-\displaystyle{\frac{c_{\alpha}}{\overline{c}_{\alpha+\eta,\eta}}}-\displaystyle{\frac{\epsilon}{2}}\right)x_{2n+1}+3\ \!x_{2n}+2+\underline{x}_{\alpha,\eta}+2\ \!x_{2n+1}^{3/4}\vspace{5pt}\\
\qquad\qquad\qquad\qquad\le\left(\displaystyle{\frac{\epsilon}{8}}-\displaystyle{\frac{\epsilon}{2}}\right)x_{2n+1}+4\ \!x_{2n}+2+\underline{x}_{\alpha,\eta}+2\ \!x_{2n+1}^{3/4}\le-\displaystyle{\frac{\epsilon}{8}}\ \!x_{2n+1}\le-A
\eaa$$
from (\ref{calphabeta}), (\ref{choiceN0}) and (\ref{tauznxn}), whence
\be\label{taux2n1}
u(\tau(x_{2n+1}),\cdot)\ge\varphi_{\alpha}(-A)-\frac{\eta}{8}\ge1-\frac{\eta}{4}\ge 1-\eta\ \hbox{ in }(-\infty,x_{2n+1}-\epsilon'x_{2n+1}-X]
\ee
from (\ref{choiceA}), (\ref{uphi}) and the monotonicity of $\varphi_{\alpha}$.

\subsubsection*{From time $t=\tau(x_{2n+1})$ to time $t=\tau(x_{2n+2})$}

{\underline{Step 1: lower bound of $\tau(x)$ for $x\ge x_{2n+1}$}}. Choose any integer $n$ such that $n\ge N_0$. There holds
$$u(\tau(x_{2n+1}),\cdot)\le\beta+\eta\ \hbox{ in }[x_{2n+1}+x_{2n+1}^{3/4},+\infty)$$
from (\ref{betaeta}) and (\ref{choiceN0}). Furthermore, $u(\tau(x_{2n+1}),\cdot)\le 1$, $\overline{\varphi}_{\beta+\eta,\eta}\ge\beta+\eta$ in $\R$, and $\overline{\varphi}_{\beta+\eta,\eta}\ge1$ in $(-\infty,-A]$ from the third assertion in (\ref{choiceA}). Thus,
$$\forall\ x\in\R,\quad u(\tau(x_{2n+1}),x)\le\overline{\varphi}_{\beta+\eta,\eta}\left(x-(x_{2n+1}+x_{2n+1}^{3/4})-A\right).$$
Since $\overline{\varphi}_{\beta+\eta,\eta}(x-\overline{c}_{\beta+\eta,\eta}t)$ is a supersolution of the equation satisfied by $u$, the maximum principle implies that
$$\forall\,(t,x)\!\in\![\tau(x_{2n+1}),+\infty)\times\R,\  u(t,x)\le\overline{\varphi}_{\beta+\eta,\eta}\!\left(\!x\!-\!x_{2n+1}\!-\!x_{2n+1}^{3/4}\!-\!A\!-\!\overline{c}_{\beta+\eta,\eta}(t\!-\!\tau(x_{2n+1}))\!\right)\!.$$
In particular, by choosing any $x\ge x_{2n+1}$ and $t=\tau(x)\ge\tau(x_{2n+1})$, one gets that
$$\theta=u(\tau(x),x)\le\overline{\varphi}_{\beta+\eta,\eta}\left(x-x_{2n+1}-x_{2n+1}^{3/4}-A-\overline{c}_{\beta+\eta,\eta}(\tau(x)-\tau(x_{2n+1}))\right)\!.$$
Since $\overline{\varphi}_{\beta+\eta,\eta}(0)=\theta$ and the function $\overline{\varphi}_{\beta+\eta,\eta}$ is decreasing, the argument of $\overline{\varphi}_{\beta+\eta,\eta}$ in the above formula is nonpositive, whence
\be\label{tauinf2}
\forall\ x\ge x_{2n+1},\quad\tau(x)\ge\frac{x-x_{2n+1}-x_{2n+1}^{3/4}-A}{\overline{c}_{\beta+\eta,\eta}}+\tau(x_{2n+1})\ge\frac{x-3\ \!x_{2n+1}}{\overline{c}_{\beta+\eta,\eta}}+\tau(x_{2n+1})
\ee
from (\ref{choiceN0}).\hfill\break

\noindent{\underline{Step 2: upper bound of $\tau(x)$ for $x\in[x_{2n+1}+x_{2n+1}^{3/4},x_{2n+2}-2\ \!x_{2n+2}^{3/4}]$}}. Let $n$ be any given integer such that $n\ge N_0$, and let $Y\ge 0$ be such that
\be\label{hyprec2}
u(\tau(x_{2n+1}),\cdot)\ge 1-\frac{\eta}{4}\ \hbox{ in }(-\infty,x_{2n+1}-\epsilon'x_{2n+1}-Y].
\ee
We are going to estimate from below, by suitable travelling fronts, the solution $u$ on the time intervals $[\tau(x_{2n+1}),\tau(x_{2n+1}+x_{2n+1}^{3/4})]$ and $[\tau(x_{2n+1}+x_{2n+1}^{3/4}),\tau(x_{2n+2}-2\ \!x_{2n+2}^{3/4})]$.\par 
Remember that $u(\tau(x_{2n+1}),\cdot)\ge\alpha$ and $\underline{\varphi}_{\alpha-\eta,\eta/4}\le1-\eta/4$ in $\R$ and that $\underline{\varphi}_{\alpha-\eta,\eta/4}\le\alpha$ in $[A,+\infty)$ from the fourth assertion in (\ref{choiceA}). Thus,
$$\forall\ x\in\R,\quad u(\tau(x_{2n+1}),x)\ge\underline{\varphi}_{\alpha-\eta,\eta/4}\left(x-(x_{2n+1}-\epsilon'x_{2n+1}-Y)+A\right).$$
Since $\underline{\varphi}_{\alpha-\eta,\eta/4}(x-\underline{c}_{\alpha-\eta,\eta/4}t)$ is a subsolution of the equation satisfied by $u$, the maximum principle implies that, for all $(t,x)\!\in\![\tau(x_{2n+1}),+\infty)\times\R$,
\be\label{varphiY}
u(t,x)\ge\underline{\varphi}_{\alpha-\eta,\eta/4}\left(x-x_{2n+1}+\epsilon'x_{2n+1}+Y+A-\underline{c}_{\alpha-\eta,\eta/4}(t-\tau(x_{2n+1}))\right)\!.
\ee\par
Let us now find a better subsolution of $u$ for times larger than $\tau(x_{2n+1}+x_{2n+1}^{3/4})$. It follows from (\ref{vimpair}) --applied at $x=x_{2n+1}$-- and the inequality $u\ge v$, that
$$u(\tau(x_{2n+1}+x_{2n+1}^{3/4}),\cdot)\ge v(\tau(x_{2n+1}+x_{2n+1}^{3/4}),\cdot)\ge\beta-\frac{\eta}{2}\ \hbox{ in }[x_{2n+1}+x_{2n+1}^{3/4},x_{2n+2}-x_{2n+2}^{3/4}].$$
Since $u(\tau(x_{2n+1}+x_{2n+1}^{3/4}),\cdot)\ge\theta\ge\beta-\eta/2$ in $(-\infty,x_{2n+1}+x_{2n+1}^{3/4}]$, one gets that
$$u(\tau(x_{2n+1}+x_{2n+1}^{3/4}),\cdot)\ge\beta-\frac{\eta}{2}\ \hbox{ in }(-\infty,x_{2n+2}-x_{2n+2}^{3/4}].$$
Furthermore, since $\tau$ is increasing and $\underline{\varphi}_{\alpha-\eta,\eta/4}$ is decreasing, it follows from the fifth assertion in~(\ref{choiceA}) and from ~(\ref{varphiY}) that
$$u(\tau(x_{2n+1}+x_{2n+1}^{3/4}),\cdot)\ge\underline{\varphi}_{\alpha-\eta,\eta/4}(-A)\ge1-\frac{\eta}{2}\ \hbox{ in }(-\infty,x_{2n+1}-\epsilon'x_{2n+1}-Y-2\ \!A].$$
Since $\underline{\varphi}_{\beta-\eta,\eta/2}\le\beta-\eta/2$ in $[A,+\infty)$ by virtue of the sixth assertion in (\ref{choiceA}), and since $\underline{\varphi}_{\beta-\eta,\eta/2}\le1-\eta/2$ in~$\R$, it resorts from the last two formulas that
$$u(\tau(x_{2n+1}+x_{2n+1}^{3/4}),\cdot)\ge\underline{\varphi}_{\beta-\eta,\eta/2}\left(\cdot-(x_{2n+1}\!-\!\epsilon'x_{2n+1}\!-\!Y\!-\!2\ \!A)\!+\!A\right)\hbox{ in }(-\infty,x_{2n+2}\!-\!x_{2n+2}^{3/4}].$$
On the other hand, there holds $x_{2n+1}\le x_{2n+2}-2\ \!x_{2n+2}^{3/4}\le x_{2n+2}-x_{2n+2}^{3/4}$ and
$$(x_{2n+2}-2\ \!x_{2n+2}^{3/4})+(x_{2n+2}-2\ \!x_{2n+2}^{3/4})^{3/4}\le x_{2n+2}-x_{2n+2}^{3/4}.$$
It follows then from (\ref{vimpair}) applied at $x=x_{2n+2}-2\ \!x_{2n+2}^{3/4}$ and from the monotonicity of $\tau$ that
$$\forall\ t\in[0,\tau(x_{2n+2}-2\ \!x_{2n+2}^{3/4})],\quad u(t,x_{2n+2}-x_{2n+2}^{3/4})\ge v(t,x_{2n+2}-x_{2n+2}^{3/4})\ge\beta-\frac{\eta}{2}.$$
Set
$$T'=\min\!\left(\!\tau(x_{2n+2}\!-\!2\ \!x_{2n+2}^{3/4}),\tau(x_{2n+1}\!+\!x_{2n+1}^{3/4})\!+\!\frac{x_{2n+2}\!-\!x_{2n+2}^{3/4}\!-\!x_{2n+1}\!+\!\epsilon'x_{2n+1}\!+\!Y\!+\!2\ \!A}{\underline{c}_{\beta-\eta,\eta/2}}\right)\!.$$
Observe that $T'\in[\tau(x_{2n+1}+x_{2n+1}^{3/4}),\tau(x_{2n+2}-2\ \!x_{2n+2}^{3/4})]$ and that
$$\underline{\varphi}_{\beta-\eta,\eta/2}\!\left(x_{2n+2}\!-\!x_{2n+2}^{3/4}\!-\!x_{2n+1}\!+\!\epsilon'x_{2n+1}\!+\!Y\!+\!3\ \!A\!-\!\underline{c}_{\beta-\eta,\eta/2}(t\!-\!\tau(x_{2n+1}+x_{2n+1}^{3/4}))\right)\!\le\beta\!-\!\frac{\eta}{2}$$
for all $t\in[\tau(x_{2n+1}+x_{2n+1}^{3/4}),T']$, because of the sixth assertion in (\ref{choiceA}) and the monotoni\-ci\-ty of $\underline{\varphi}_{\beta-\eta,\eta/2}$. Eventually, since $\underline{\varphi}_{\beta-\eta,\eta/2}(x-\underline{c}_{\beta-\eta,\eta/2}t)$ is a subsolution of the equation satisfied by $u$, the maximum principle applied in $[\tau(x_{2n+1}+x_{2n+1}^{3/4}),T']\times(-\infty,x_{2n+2}-x_{2n+2}^{3/4}]$ yields
\be\label{maxprin}
u(t,x)\ge\underline{\varphi}_{\beta-\eta,\eta/2}\!\left(x-\!x_{2n+1}\!+\!\epsilon'x_{2n+1}\!+\!Y\!+\!3\ \!A\!-\!\underline{c}_{\beta-\eta,\eta/2}(t\!-\!\tau(x_{2n+1}+x_{2n+1}^{3/4}))\right)
\ee
for all $(t,x)\in[\tau(x_{2n+1}+x_{2n+1}^{3/4}),T']\times(-\infty,x_{2n+2}-x_{2n+2}^{3/4}]$.\par
Pick any $x\in[x_{2n+1}+x_{2n+1}^{3/4},x_{2n+2}-2\ \!x_{2n+2}^{3/4}]$, set
$$t(x)=\tau(x_{2n+1}+x_{2n+1}^{3/4})+\frac{x-x_{2n+1}+\epsilon'x_{2n+1}+Y+3\ \!A}{\underline{c}_{\beta-\eta,\eta/2}}$$
and assume that $\tau(x)>t(x)$. Then $t(x)\le\tau(x)\le\tau(x_{2n+2}-2\ \!x_{2n+2}^{3/4})$ since $\tau$ is increasing. On the other hand,
\be\label{tx}\baa{rcl}
t(x) & \le & \tau(x_{2n+1}\!+\!x_{2n+1}^{3/4})\!+\!\displaystyle{\frac{x_{2n+2}\!-\!2\ \!x_{2n+2}^{3/4}\!-\!x_{2n+1}\!+\!\epsilon'x_{2n+1}\!+\!Y\!+\!3\ \!A}{\underline{c}_{\beta-\eta,\eta/2}}}\vspace{5pt}\\
& \le & \tau(x_{2n+1}\!+\!x_{2n+1}^{3/4})\!+\!\displaystyle{\frac{x_{2n+2}\!-\!x_{2n+2}^{3/4}\!-\!x_{2n+1}\!+\!\epsilon'x_{2n+1}\!+\!Y\!+\!2\ \!A}{\underline{c}_{\beta-\eta,\eta/2}}}\eaa
\ee
since $-x_{2n+2}^{3/4}+A\le 0$, because of (\ref{choiceN0}). Thus, $t(x)\le T'$. Observe also that
$$t(x)\ge\tau(x_{2n+1}+x_{2n+1}^{3/4})$$
by definition of $t(x)$ and since $x\ge x_{2n+1}+x_{2n+1}^{3/4}\ge x_{2n+1}$ and all parameters~$Y$ and $A$ are nonnegative. One can then apply (\ref{maxprin}) at the point $(t(x),x)$ and one gets
$$u(t(x),x)\ge\underline{\varphi}_{\beta-\eta,\eta/2}(0)=\theta,$$
whence $\tau(x)\le t(x)$, owing to the definition of $\tau(x)$. As a consequence, the assumption $\tau(x)>t(x)$ cannot hold and one concludes that
\be\label{tausup2}\baa{l}
\forall\ x\in[x_{2n+1}+x_{2n+1}^{3/4},x_{2n+2}-2\ \!x_{2n+2}^{3/4}],\vspace{5pt}\\
\qquad\qquad\baa{rcl} \tau(x) & \le & t(x)=\tau(x_{2n+1}+x_{2n+1}^{3/4})+\displaystyle{\frac{x-x_{2n+1}+\epsilon'x_{2n+1}+Y+3\ \!A}{\underline{c}_{\beta-\eta,\eta/2}}}\vspace{5pt}\\
& \le & \tau(x_{2n+1}+x_{2n+1}^{3/4})+\displaystyle{\frac{x+Y+3\ \!A}{\underline{c}_{\beta-\eta,\eta/2}}}.\eaa\eaa
\ee

\noindent{\underline{Step 3: estimate of $\tau(z_{2n+1})$ and lower bound of $u$ on the left of $z_{2n+1}$}}. It follows from (\ref{taux}), (\ref{tauinf2}), (\ref{tausup2}) and the inequality
$$x_{2n+1}+x_{2n+1}^{3/4}\le z_{2n+1}\le x_{2n+2}-2\ \!x_{2n+2}^{3/4}$$
that
\be\label{tauz2n1}\baa{l}
\displaystyle{\frac{z_{2n+1}-3\ \!x_{2n+1}}{\overline{c}_{\beta+\eta,\eta}}}\le\displaystyle{\frac{z_{2n+1}-3\ \!x_{2n+1}}{\overline{c}_{\beta+\eta,\eta}}}+\tau(x_{2n+1})\le\cdots\vspace{5pt}\\
\qquad\qquad\qquad\qquad\cdots\le\tau(z_{2n+1})\le\displaystyle{\frac{x_{2n+1}+x_{2n+1}^{3/4}+M}{c_{\alpha}}}+\displaystyle{\frac{z_{2n+1}+Y+3\ \!A}{\underline{c}_{\beta-\eta,\eta/2}}},\eaa
\ee
provided that (\ref{hyprec2}) holds. Furthermore,
$$\baa{l}
\tau(x_{2n+1}+x_{2n+1}^{3/4})\le\tau(z_{2n+1})\le\tau(x_{2n+2}-2\ \!x_{2n+2}^{3/4})\le t(x_{2n+2}-2\ \!x_{2n+2}^{3/4})\le\cdots\vspace{5pt}\\
\qquad\qquad\qquad\qquad\cdots\le\tau(x_{2n+1}\!+\!x_{2n+1}^{3/4})\!+\!\displaystyle{\frac{x_{2n+2}\!-\!x_{2n+2}^{3/4}\!-\!x_{2n+1}\!+\!\epsilon'x_{2n+1}\!+\!Y\!+\!2\ \!A}{\underline{c}_{\beta-\eta,\eta/2}}}\eaa$$
from (\ref{tx}) and (\ref{tausup2}). Thus, $\tau(z_{2n+1})\le T'$. The inequality (\ref{maxprin}) and the monotonicity of $\underline{\varphi}_{\beta-\eta,\eta/2}$ then imply that, for all $x\in(-\infty,z_{2n+1}-\epsilon'z_{2n+1}-Y]\ (\subset(-\infty,x_{2n+2}-x_{2n+2}^{3/4}])$,
$$\baa{l}
u(\tau(z_{2n+1}),x)\ge\underline{\varphi}_{\beta-\eta,\eta/2}(z_{2n+1}-\epsilon'z_{2n+1}-x_{2n+1}+\epsilon'x_{2n+1}+3\ \!A\cdots\vspace{5pt}\\
\qquad\qquad\qquad\qquad\qquad\qquad\qquad\qquad\cdots-\underline{c}_{\beta-\eta,\eta/2}(\tau(z_{2n+1})-\tau(x_{2n+1}+x_{2n+1}^{3/4}))).\eaa$$
But it follows from (\ref{taux}), (\ref{calphabeta}), (\ref{choiceN0}) and (\ref{tauz2n1}) that
$$\baa{l}
z_{2n+1}-\epsilon'z_{2n+1}-x_{2n+1}+\epsilon'x_{2n+1}+3\ \!A-\underline{c}_{\beta-\eta,\eta/2}(\tau(z_{2n+1})-\tau(x_{2n+1}+x_{2n+1}^{3/4}))\vspace{5pt}\\
\qquad\qquad\le\left(1-\displaystyle{\frac{\underline{c}_{\beta-\eta,\eta/2}}{\overline{c}_{\beta+\eta,\eta}}}-\displaystyle{\frac{\epsilon}{2}}\right)z_{2n+1}+3\ \!A+\displaystyle{\frac{3\ \!x_{2n+1}\ \!\underline{c}_{\beta-\eta,\eta/2}}{\overline{c}_{\beta+\eta,\eta}}}+\displaystyle{\frac{c_{\beta}}{c_{\alpha}}}(x_{2n+1}+x_{2n+1}^{3/4}+M)\vspace{5pt}\\
\qquad\qquad\le\left(\displaystyle{\frac{\epsilon}{8}}-\displaystyle{\frac{\epsilon}{2}}\right)z_{2n+1}+\left(6+\displaystyle{\frac{3\ \!c_{\beta}}{c_{\alpha}}}\right)x_{2n+1}\le-\displaystyle{\frac{\epsilon}{8}}z_{2n+1}\le-A,\eaa$$
whence
\be\label{tauz2n1bis}
u(\tau(z_{2n+1}),\cdot)\ge\underline{\varphi}_{\beta-\eta,\eta/2}(-A)\ge1-\frac{3\eta}{4}\ge1-\eta\ \hbox{ in }(-\infty,z_{2n+1}-\epsilon'z_{2n+1}-Y].
\ee
from the seventh assertion in (\ref{choiceA}).\hfill\break

\noindent{\underline{Step 4: estimate of $u$ on the left of $x_{2n+2}$ at time $\tau(x_{2n+2}-2\ \!x_{2n+2}^{3/4})$}}. With similar arguments as above, one has 
$$\tau(x_{2n+1}+x_{2n+2}^{3/4})\le\tau(x_{2n+2}-2\ \!x_{2n+2}^{3/4})\le T'$$
and
$$x_{2n+2}-\rho\epsilon x_{2n+2}-Y\le x_{2n+2}-x_{2n+2}^{3/4},$$
since $Y\ge 0$ in (\ref{hyprec2}) and since $\rho\,\epsilon\,x_{2n+2}^{1/4}\ge 1$ from (\ref{choiceN0}). Thus, inequality (\ref{maxprin}) and the monotonicity of $\underline{\varphi}_{\beta-\eta,\eta/2}$ imply that, for all $x\in(-\infty,x_{2n+2}-\rho\epsilon x_{2n+2}-Y]$,
$$\baa{l}
u(\tau(x_{2n+2}-2\ \!x_{2n+2}^{3/4}),x)\ge\underline{\varphi}_{\beta-\eta,\eta/2}\Big(x_{2n+2}-\rho\epsilon x_{2n+2}-x_{2n+1}+\epsilon'x_{2n+1}+3\ \!A\cdots\vspace{5pt}\\
\qquad\qquad\qquad\qquad\qquad\qquad\qquad\qquad\cdots-\underline{c}_{\beta-\eta,\eta/2}(\tau(x_{2n+2}\!-\!2\ \!x_{2n+2}^{3/4})\!-\!\tau(x_{2n+1}\!+\!x_{2n+1}^{3/4}))\Big).\eaa$$
But, as in Step~3, it follows from (\ref{taux}), (\ref{calphabeta}), (\ref{choiceN0}) and (\ref{tauinf2}) applied at $x_{2n+2}-2\ \!x_{2n+2}^{3/4}$, that
$$\baa{l}
x_{2n+2}\!-\!\rho\epsilon x_{2n+2}\!-\!x_{2n+1}\!+\!\epsilon'x_{2n+1}\!+\!3\ \!A\!-\!\underline{c}_{\beta-\eta,\eta/2}(\tau(x_{2n+2}\!-\!2\ \!x_{2n+2}^{3/4})-\tau(x_{2n+1}+x_{2n+1}^{3/4}))\vspace{5pt}\\
\qquad\qquad\le\left(\displaystyle{\frac{\rho\epsilon}{4}}+\displaystyle{\frac{\rho\epsilon}{2}}-\rho\epsilon\right)x_{2n+2}+\displaystyle{\frac{2\ \!\underline{c}_{\beta-\eta,\eta/2}\ \!x_{2n+2}^{3/4}}{\overline{c}_{\beta+\eta,\eta}}}\le-\displaystyle{\frac{\rho\epsilon}{4}}x_{2n+2}+2\ \!x_{2n+2}^{3/4}\le-A,\eaa$$
whence
\be\label{taux2n20}
u(\tau(x_{2n+2}-2\ \!x_{2n+2}^{3/4}),\cdot)\ge\underline{\varphi}_{\beta-\eta,\eta/2}(-A)\ge1-\frac{3\eta}{4}\ \hbox{ in }\left(-\infty,x_{2n+2}-\rho\epsilon x_{2n+2}-Y\right].
\ee

\noindent{\underline{Step 5: estimates of $\tau(x_{2n+2})$ and of $u$ on the left of $x_{2n+2}$ at time $\tau(x_{2n+2})$}}. Remember that $u(\tau(x_{2n+2}-2\ \!x_{2n+2}^{3/4}),\cdot)\ge\alpha$ in $\R$, that $\underline{\varphi}_{\alpha-\eta,3\eta/4}\le1-3\eta/4$ in $\R$ and that $\underline{\varphi}_{\alpha-\eta,3\eta/4}\le\alpha$ in $[A,+\infty)$ from the eighth assertion in (\ref{choiceA}). Thus,
$$\forall\ x\in\R,\quad u(\tau(x_{2n+2}-2\ \!x_{2n+2}^{3/4}),x)\ge\underline{\varphi}_{\alpha-\eta,3\eta/4}\left(x-(x_{2n+2}-\rho\epsilon x_{2n+2}-Y)+A\right)$$
from (\ref{taux2n20}), provided that (\ref{hyprec2}) holds. Since $\underline{\varphi}_{\alpha-\eta,3\eta/4}(x-\underline{c}_{\alpha-\eta,3\eta/4}t)$ is a subsolution of the equation satisfied by $u$, the maximum principle yields
\be\label{maxprin2}
u(t,x)\ge\underline{\varphi}_{\alpha-\eta,3\eta/4}\!\left(\!x\!-\!x_{2n+2}\!+\!\rho\epsilon x_{2n+2}\!+\!Y\!+\!A\!-\!\underline{c}_{\alpha-\eta,3\eta/4}(t\!-\!\tau(x_{2n+2}\!-\!2\ \!x_{2n+2}^{3/4}))\!\right)
\ee
for all $(t,x)\in[\tau(x_{2n+2}-2\ \!x_{2n+2}^{3/4}),+\infty)\times\R$. In particular, at
$$t=\tau(x_{2n+2})\ge\tau(x_{2n+2}-2\ \!x_{2n+2}^{3/4})\ \hbox{ and }\ x=x_{2n+2},$$
one gets that
$$\theta\!=\!u(\tau(x_{2n+2}),x_{2n+2})\!\ge\!\underline{\varphi}_{\alpha-\eta,3\eta/4}\!\!\left(\!\rho\epsilon x_{2n+2}\!+\!\!Y\!\!+\!\!A\!-\!\underline{c}_{\alpha-\eta,3\eta/4}(\tau(x_{2n+2})\!-\!\tau(x_{2n+2}\!-\!2\ \!x_{2n+2}^{3/4}))\!\right)\!\!,$$
whence
$$\tau(x_{2n+2})\le\tau(x_{2n+2}-2\ \!x_{2n+2}^{3/4})+\frac{\rho\epsilon x_{2n+2}+Y+A}{\underline{c}_{\alpha-\eta,3\eta/4}}$$
since $\underline{\varphi}_{\alpha-\eta,3\eta/4}$ is decreasing and $\underline{\varphi}_{\alpha-\eta,3\eta/4}(0)=\theta$. It follows then from (\ref{taux}), (\ref{tauinf2}) and~(\ref{tausup2}) applied at $x=x_{2n+2}-2\ \!x_{2n+2}^{3/4}$, that
\be\label{taux2n2}\baa{l}
\displaystyle{\frac{x_{2n+2}-3\ \!x_{2n+1}}{\overline{c}_{\beta+\eta,\eta}}}\le\displaystyle{\frac{x_{2n+2}-3\ \!x_{2n+1}}{\overline{c}_{\beta+\eta,\eta}}}+\tau(x_{2n+1})\le\tau(x_{2n+2})\le\cdots\vspace{5pt}\\
\qquad\cdots\le\displaystyle{\frac{x_{2n+1}+x_{2n+1}^{3/4}+M}{c_{\alpha}}}+\displaystyle{\frac{x_{2n+2}-2\ \!x_{2n+2}^{3/4}+Y+3\ \!A}{\underline{c}_{\beta-\eta,\eta/2}}}+\displaystyle{\frac{\rho\epsilon x_{2n+2}+Y+A}{\underline{c}_{\alpha-\eta,3\eta/4}}}.\eaa
\ee\par
Lastly, inequality (\ref{maxprin2}) applied at $t=\tau(x_{2n+2})\ge\tau(x_{2n+2}-2\ \!x_{2n+2}^{3/4})$ implies that
\be\label{taux2n2bis}\baa{rcl}
\forall\ x\le x_{2n+2}\!-\!\epsilon'x_{2n+2}\!-\!Y,\ \ \ u(\tau(x_{2n+2}),x) & \!\!\!\ge\!\!\! & \underline{\varphi}_{\alpha-\eta,3\eta/4}\!\left(\!-\displaystyle{\frac{\epsilon}{2}x_{2n+2}}\!+\!\rho\epsilon x_{2n+2}\!+\!A\!\right)\vspace{5pt}\\
& \!\!\!\ge\!\!\! & \underline{\varphi}_{\alpha-\eta,3\eta/4}(-A)\ge1-\eta\eaa
\ee
since $\underline{\varphi}_{\alpha-\eta,3\eta/4}$ is decreasing and because of the last assertion in (\ref{choiceA}) and because of~(\ref{choiceN0}).

\subsubsection*{Conclusion of the proof of Lemma~\ref{speeds}}

As already underlined, for all $N\ge N_0$, the estimates
$$\left\{\baa{ll}
u(\tau(x_{2n},\cdot)\le\alpha+\eta &\hbox{ in }[x_{2n}+\epsilon x_{2n},x_{2n+1}-\epsilon x_{2n+1}],\vspace{5pt}\\
u(\tau(z_{2n},\cdot)\le\alpha+\eta & \hbox{ in }[z_{2n}+\epsilon x_{2n},x_{2n+1}-\epsilon x_{2n+1}],\vspace{5pt}\\
|u(\tau(x_{2n+1},\cdot)-\beta|\le\eta &\hbox{ in }[x_{2n+1}+\epsilon x_{2n+1},x_{2n+2}-\epsilon x_{2n+2}],\vspace{5pt}\\
|u(\tau(z_{2n+1},\cdot)-\beta|\le\eta & \hbox{ in }[z_{2n+1}+\epsilon z_{2n+1},x_{2n+2}-\epsilon x_{2n+2}]\eaa\right.$$
follow from (\ref{choiceN0}), (\ref{alphaeta}) and (\ref{right}).\par
Now, since $u(\tau(x_{2N_0}),-\infty)=1$, there exists a nonnegative real number $X_{2N_0}=X_{2N_0}(\epsilon,\eta)\ge 0$ such that
$$u(\tau(x_{2N_0}),\cdot)\ge 1-\eta\ \hbox{ in }(-\infty,x_{2N_0}-\epsilon'x_{2N_0}-X_{2N_0}].$$
In other words, $X_{2N_0}$ plays the role of $X$ in (\ref{hyprec}), with $n=N_0$. It follows then from~(\ref{tauz2n}) and~(\ref{taux2n1}) that
\be\label{left1}
u(\tau(z_{2N_0}),\cdot)\ge 1-\frac{\eta}{4}\ge 1-\eta\ \hbox{ in }(-\infty,z_{2N_0}-\epsilon'z_{2N_0}-X_{2N_0}]
\ee
and that
\be\label{left2}
u(\tau(x_{2N_0+1},\cdot)\ge 1-\frac{\eta}{4}\ge 1-\eta\ \hbox{ in }(-\infty,x_{2N_0+1}-\epsilon'x_{2N_0+1}-X_{2N_0}].
\ee
As a consequence, property (\ref{hyprec2}) is fulfilled with $Y=X_{2N_0}$ and $n=N_0$. It follows then from (\ref{tauz2n1bis}) and (\ref{taux2n2bis}) that
\be\label{left3}
u(\tau(z_{2N_0+1}),\cdot)\ge 1-\frac{3\eta}{4}\ge 1-\eta\ \hbox{ in }(-\infty,z_{2N_0+1}-\epsilon'z_{2N_0+1}-X_{2N_0}]
\ee
and that
\be\label{left4}
u(\tau(x_{2N_0+2},\cdot)\ge 1-\eta\ \hbox{ in }(-\infty,x_{2N_0+2}-\epsilon'x_{2N_0+2}-X_{2N_0}].
\ee
By an immediate induction, one gets that the above four estimates (\ref{left1})-(\ref{left4}) hold for all $n\ge N_0$. Hence, since $\epsilon'=\epsilon/2>0$ and $\lim_{m\to+\infty}x_m=\lim_{m\to+\infty}z_m=+\infty$, there exists an integer $N_1=N_1(\epsilon,\eta)\ge N_0$ such that
$$\forall\ n\ge N_1,\left\{\baa{l}
u(\tau(x_{2n}),\cdot)\ge 1-\eta\ \hbox{ in }(-\infty,x_{2n}-\epsilon x_{2n}],\vspace{5pt}\\
u(\tau(z_{2n}),\cdot)\ge 1-\eta\ \hbox{ in }(-\infty,z_{2n}-\epsilon z_{2n}],\vspace{5pt}\\
u(\tau(x_{2n+1}),\cdot)\ge 1-\eta\ \hbox{ in }(-\infty,x_{2n+1}-\epsilon x_{2n+1}],\vspace{5pt}\\
u(\tau(z_{2n+1}),\cdot)\ge 1-\eta\ \hbox{ in }(-\infty,z_{2n+1}-\epsilon z_{2n+1}].\eaa\right.$$
Furthermore, by an immediate induction, it also follows that the estimates (\ref{tauznxn}), (\ref{tauz2n1}) and (\ref{taux2n2}) hold for all $n\ge N_0$ with $X=Y=X_{2N_0}$. Therefore, since $\lim_{m\to+\infty}x_m=\lim_{m\to+\infty}x_{m+1}/x_m=\lim_{m\to+\infty}z_m/x_m=+\infty$, one gets that
$$\left\{\baa{l}
c_{\alpha}\le\displaystyle{\mathop{\liminf}_{n\to+\infty}}\ \displaystyle{\frac{z_{2n}}{\tau(z_{2n})}}\le\displaystyle{\mathop{\limsup}_{n\to+\infty}}\ \displaystyle{\frac{z_{2n}}{\tau(z_{2n})}}\le\overline{c}_{\alpha+\eta,\eta},\vspace{5pt}\\
c_{\alpha}\le\displaystyle{\mathop{\liminf}_{n\to+\infty}}\ \displaystyle{\frac{x_{2n+1}}{\tau(x_{2n+1})}}\le\displaystyle{\mathop{\limsup}_{n\to+\infty}}\ \displaystyle{\frac{x_{2n+1}}{\tau(x_{2n+1})}}\le\overline{c}_{\alpha+\eta,\eta},\vspace{5pt}\\
\underline{c}_{\beta-\eta,\eta/2}\le\displaystyle{\mathop{\liminf}_{n\to+\infty}}\ \displaystyle{\frac{z_{2n+1}}{\tau(z_{2n+1})}}\le\displaystyle{\mathop{\limsup}_{n\to+\infty}}\ \displaystyle{\frac{z_{2n+1}}{\tau(z_{2n+1})}}\le\overline{c}_{\beta+\eta,\eta},\vspace{5pt}\\
\left(\underline{c}_{\beta-\eta,\eta/2}^{-1}+\rho\ \!\epsilon\ \!\underline{c}_{\alpha-\eta,3\eta/4}^{-1}\right)^{-1}\le\displaystyle{\mathop{\liminf}_{n\to+\infty}}\ \displaystyle{\frac{x_{2n+2}}{\tau(x_{2n+2})}}\le\displaystyle{\mathop{\limsup}_{n\to+\infty}}\ \displaystyle{\frac{x_{2n+2}}{\tau(x_{2n+2})}}\le\overline{c}_{\beta+\eta,\eta}.\eaa\right.$$
Because of (\ref{calphabeta}), there exists an integer $N=N(\epsilon,\eta)\ge N_1$ such that
$$\forall\ n\ge N,\quad
\left|\displaystyle{\frac{x_{2n}}{\tau(x_{2n})}}-c_{\beta}\right|+\left|\displaystyle{\frac{z_{2n}}{\tau(z_{2n})}}-c_{\alpha}\right|+\left|\displaystyle{\frac{x_{2n+1}}{\tau(x_{2n+1})}}-c_{\alpha}\right|+\left|\displaystyle{\frac{z_{2n+1}}{\tau(z_{2n+1})}}-c_{\beta}\right|\le\epsilon.$$
That completes the proof of Lemma~\ref{speeds}.\hfill$\Box$

\begin{rem}\label{remNL}{\rm The behaviour of the solution $u$ in the region where it is less than $\theta$ is close in some sense to that of the solution $v$ of the heat equation, as expected. The function $v$ oscillates at large time between $\alpha$ and $\beta$, infinitely many times: such a nontrivial dynamics is well-known for the heat equation, see \cite{ce}. However, the difficulty in the above proof came from the nonlinear reaction term $f(u)$ and from the estimates of the position and average speed of the solution $u$ as time runs.}
\end{rem}

%%%%%%%%%%%%%%%%%%%%%%%%%%%%%%%%%%%%%%%%
%%%%%%%%%%%%%%%%%%%%%%%%%%%%%%%%%%%%%%%%

\end{document}